\documentclass{elsarticle}
\usepackage{mathtools}
\usepackage{enumitem}
\usepackage{xcolor}
\usepackage{amsmath}
\usepackage{amsthm}
\usepackage{upgreek}
\usepackage{graphics}
\usepackage{makeidx}
\usepackage{color}
\usepackage{graphicx}
\usepackage{amssymb}
\usepackage{float}
\usepackage[english]{babel}
\usepackage{caption}
\usepackage{subcaption}
\usepackage{longtable}
\usepackage[utf8]{inputenc}
\usepackage[labelfont=bf]{caption}
\usepackage[singlelinecheck=off]{caption}

\usepackage{amsthm}

\theoremstyle{definition}

\makeindex
\makeglossary
\begin{document}

\begin{frontmatter}
\title{Formation of stationary periodic patterns in a model of two competing populations with chemotaxis}

	\author{Valentina Bucur and Bakhtier Vasiev\corref{email}}

\cortext[email]{bnvasiev@liverpool.ac.uk}

\address{Department of Mathematical Sciences, Mathematical Sciences Building, University of Liverpool, Liverpool L69 7ZL, UK}

\begin{abstract}
	One of the classical models in mathematical biology is the Lotka-Volterra competition model, describing the dynamics of two populations competing for resources. Two possible regimes in this system are given by their coexistence or extinction of a weaker population. In a distributed system with diffusive spatial coupling, travelling fronts occur, corresponding to transitions between stationary states. In this work we will consider the competition model extended by extra interaction between involved populations which is given by chemotactic coupling, namely, assuming that one species produces a chemical agent which causes the taxis of another species. It is known that in a one-species model (i.e. Keller-Segel model) production of chemoattractor results in formation of stationary periodic (or Turing-type) patterns. In this work, we will investigate conditions for the formation of stationary periodic patterns in a two-species competition model with chemotaxis. We show that in this system periodic patterns can emerge in the course of Turing-type instability (classical way) or from a stable steady state, corresponding to the extinction of one of the species, due to a finite, or over-threshold, amplitude disturbance. We study the characteristics of emerging periodic pattern, such as its amplitude and wavelength, by means of Fourier analysis. We also perform computational simulations to verify our analytical results. 
\end{abstract}
\end{frontmatter}

\section{Introduction}

Ecological systems, such as bacterial biofilms, are often composed of multiple biological species whose interactions play a crucial role in the formation and properties of these systems \cite{JMLang}. The interactions between biological populations can take various forms, including predator-prey, competition, and symbiosis \cite{Murray1}. In this work, we consider the case of competitive interactions between two populations, which is commonly modelled by the classical Lotka-Volterra model \cite{Lotka, Volterra}, represented by the following (nondimensional) equations:
\begin{equation} \label{LW}
\displaystyle\begin{cases}
\displaystyle u_t=u(1-u-b_1v),\\
\displaystyle v_t=rv(1-v-b_2u),
\end{cases}
\end{equation}
in which variables  $u$ and $v$ represent sizes of two species, parameter $r$ gives the ratio of reproduction rates of $v$ and  $u$, parameter $b_1$ represents the competition strength of $v$ on  $u$ and $b_2$ - of  $u$ on $v$. The case when $0<b_1<1$ and $0<b_2<1$ is referred to as a weak competition regime; the case of $b_1>1$ and $b_2>1$ corresponds to strong competition regime, and intermediate cases as string-weak or weak-strong regimes. It is well known that this model has four steady states and their stability is defined by model parameters as follows \cite{Hsu, Murray1, Shigesada1}:
\begin{itemize}
\item $(u_1,v_1)=(0, 0)$ (both species are extinct) is always unstable;
\item $\displaystyle (u_2,v_2)= \left(\frac{b_1-1}{b_1b_2-1}, \frac{b_2-1}{b_1b_2-1}\right)$ (coexistence) is stable if $b_1, b_2<1$ (weak competition) and unstable otherwise.
\item $(u_3,v_3)=(1, 0)$ (the second species is extinct) is stable if $b_2>1$, and unstable otherwise;
\item $(u_4,v_4)=(0, 1)$ (the first species is extinct) is stable if $b_1>1$ and unstable otherwise;
\end{itemize}
The system \eqref{LW} is often extended by adding spatial variable so that the variables $u=u(x,t)$ and $v=v(x,t)$ represent densities of species, which depend on the coordinate $x$ as well as on time $t$. Spatial coupling is commonly given by diffusion and cross-diffusion (or cross-taxis) and mathematically represented by Shigesada-Kawasaki-Teramoto, or SKT system \cite{Shigesada1}, which can be reduced to Potts-Petrovskii model \cite{Potts}: 
\begin{equation}\label{SKT}
	\begin{cases}
		\displaystyle u_t=u_{xx}-\chi_1(uv_x)_x+u(1-u-b_1v),\\
		\displaystyle v_t=dv_{xx}-\chi_2(vu_x)_x+rv(1-v-b_2u).\\
	\end{cases}
\end{equation}
Here parameter $d$ represents the rate of random motion of species $v$ relative to that of  $u$ and parameters $\chi_1$ and $\chi_2$ give the strengths of cross-taxis between species. The case of $\chi_1=\chi_2=0$ is referred to as a classical competition-diffusion model known by travelling front solutions corresponding to transitions between stationary states of system \eqref{LW}, which were studied in many papers  \cite{Hosono2, MRodrigo, Alzahrani1, Alzahrani2}. Cross-taxis terms affect the speed, direction and shape of travelling fronts \cite{WangZ} and, under certain conditions, also cause the loss of stability of equilibrium states resulting to formation of stationary periodic patterns \cite{Qing, Breden}. 

Going one step further, mathematical models have been introduced for two competing species with extra interaction mediated by chemical produced by one or both species. This chemical can affect the kinetic terms resulting, for example, to inhibition of one species by the other \cite{Dean}. It can also act as a chemotactic agent, as it is assumed in a class of so called two-species chemotaxis model, which can be represented by the following (non-dimensional) system \cite{Kelly}: 
\begin{equation}\label{twosp_mutrep}
\begin{cases}
\displaystyle u_t=D_1u_{xx}-\chi_1(uc_x)_x+u(1-u-b_1v),\\
\displaystyle v_t=D_2v_{xx}-\chi_2(vc_x)_x+rv(1-v-b_2u),\\
c_t=c_{xx}+\mu_1 u+\mu_2 v-\gamma c,
\end{cases}
\end{equation}
where parameters $D_1$ and $D_2$ define the rate of random motion of species  $u$ and $v$, $\chi_1$ and $\chi_2$ the strength of their chemotactic response to the chemical $c$, which is produced by both species (with production rates $\mu_1$ and $\mu_2$) and degrades with rate $\gamma$ \cite{HKurt}. This system models spatio-temporal evolution of two competitive species, which migrate along (chemoattraction, $\chi>0$) or against (chemorepulsion, $\chi<0$) the concentration gradient of the chemical produced by themselves. 

The impact of chemotactic terms on the shape and speed of travelling fronts, which represent transitions between the steady states of \eqref{LW}, has been in the focus of many studies of system \eqref{twosp_mutrep} and other similar systems \cite{Issa, LDongXH, TCLin, Lin, CStinner}. An analytical study of system \eqref{twosp_mutrep}, conducted in \cite{Xpan}, demonstrated that in the case of chemoattraction ($\chi_1 > 0$, $\chi_2 > 0$) and strong competition ($b_1 > 1$, $b_2 > 1$), the stable steady states can lose their stability when the chemoattraction is sufficiently strong. The loss of steady-state stability due to special coupling (given by the chemotactic term in the given case) is commonly referred to as Turing instability, with the stationary periodic patterns forming as a result known as Turing patterns \cite{Turing}. In \cite{Xpan}, it was shown that a stationary periodic pattern may arise when the steady state $(u, \mbox{ } v, \mbox{ } c) = (1, \mbox{ } 0,\mbox{ } 1)$ becomes unstable if $\chi_1$ is sufficiently large, and similarly, such patterns may arise from $(u, \mbox{ } v, \mbox{ } c) = (0, \mbox{ } 1,\mbox{ } 1)$ when $\chi_2$ is sufficiently large. Moreover, it was numerically demonstrated that in the case of weak competition, stationary periodic patterns can arise from the coexistence state if either of the parameters $\chi_1$ or $\chi_2$, defining the chemotactic response, is sufficiently large. These findings are not surprising, as by setting $u = 0$ or $v = 0$, the system \eqref{twosp_mutrep} is reduced to a one-species model, aligning with the classical result \cite{Keller, Bucur} that strong chemoattraction in a single species with logistic growth leads to Turing-type instability.

The formation of stationary periodic patterns in the case of mutual repulsion between species was reported in \cite{Li_2022}. The model presented in that paper included two chemicals, each produced by one of the bacterial species and acting as a chemorepellent to the other. It was shown that, in the case of weak competition, the coexistence state can become unstable, leading to the formation of periodic patterns when the strengths of chemorepulsion are sufficiently high. The study of these patterns in \cite{Li_2022} included an estimation of their amplitude using the amplitude equation, as well as their numerical reproduction in one- and two-dimensional domains. The phenomenon reported in \cite{Li_2022} is quite novel, as it is not observed in single-species systems, where chemorepulsion does not result to formation of periodic patterns.

In this paper, we focus on the formation of periodic patterns in system \eqref{twosp_mutrep} in the most basic scenario, where only one species produces a chemical agent which repels the second species. In line with the results reported in \cite{Li_2022}, we found that Turing-type instability occurs only in the case of weak competition, where strong chemorepulsion causes the coexistence state to become unstable, resulting in the formation of periodic patterns. We have also found that both steady states, corresponding to the extinction of one of the species (in the case of strong competition), remain stable in the presence of chemotaxis. However, we have shown that, even in this case, stationary periodic patterns can be initiated by finite perturbations. We have also analysed the wavelengths and amplitude of the forming patterns using Fourier series analysis, as done in \cite{Bucur}.

\section{Model and linear stability analysis of its steady states}

We consider a two-species competition model with chemotaxis where one of species produces chemotaxtic agent for the other. This is represented by a slightly modified version of the system \eqref{twosp_mutrep} given as:  
\begin{equation}\displaystyle\label{twospecieschem}
\displaystyle\begin{cases}
\displaystyle u_t=D_1u_{xx}-\chi(uc_x)_x+r_1u(1-u-b_1v),\\
\displaystyle v_t=D_2v_{xx}-r_2v(1-v-b_2u),\\
\displaystyle c_t=c_{xx}+v-c,
\end{cases}
\end{equation}
where the chemical agent $c$ is produced by species $v$ and affects chemotactically $u$. This system is considered on a one-dimensional domain $x\in (0,L)$ under no-flux boundary conditions. In the well-mixed case, when  the solution is homogeneous, the first two equations of this system transform into the system \eqref{LW}, last equation gives $c=v$, and therefore it has four steady states:
\[\displaystyle (u^*,\mbox{ }v^*,\mbox{ }c^*)=\left\{(0, 0, 0), \mbox{ }\left(\frac{b_1-1}{b_1b_2-1}, \frac{b_2-1}{b_1b_2-1}, \frac{b_2-1}{b_1b_2-1}\right), \mbox{ } (1, 0, 0), \mbox{ }(0, 1, 1) \right\}, \]
with their stability given by the same conditions as for the Lotka-Volterra model \eqref{LW}. 

Classical Turing instability analysis refers to the formation of stationary periodic patterns in a reaction-diffusion-advection system, when the steady states loose their stability due diffusion and advection terms \cite{Turing, Keller}. Linear stability analysis involves consideration of the perturbation given by cosinusoidal function, $\sim \exp(\lambda t)\cos(kx)$ \cite{Murray2}, so that the stability of the steady state, $(u^*,v^*,c^*)$, is defined by the signs of eigenvalues of the characteristic matrix at this state, which, for the system \eqref{twospecieschem}, is given as following:    
\begin{equation}\label{Mrda} \displaystyle M=
\begin{pmatrix}
\displaystyle -D_1k^2+r_1(1-2u^*-b_1v^*) & -r_1b_1u^* & -\chi u^*k^2 \\
\displaystyle -r_2b_2v^* & -D_2k^2+r_2(1-2v^*-b_2u^*) & 0 \\
\displaystyle 0 & 1& -k^2-1
\end{pmatrix}.
\end{equation}
The steady state is unstable if at least one of eigenvalues is positive. Since the steady state $(u_1,v_1,c_1)=(0, 0, 0)$ is always unstable and can't loose its stability due to diffusion and advection terms we don't need to consider it.  However, stationary periodic patterns can emerge from other steady states when they are stable in the well-mixed system and become unstable in the presence of diffusion and advection.

\subsection{Analysis of the coexistence state}

\begin{figure}[h]
\centering
\begin{subfigure}{.32\textwidth}
\caption{}
\includegraphics[scale=.29]{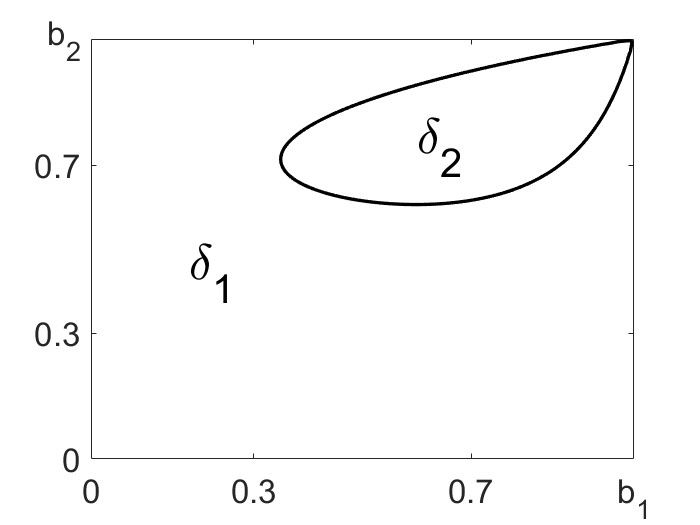}
\end{subfigure}
\hfill
\begin{subfigure}{.47\textwidth}
\caption{}
\includegraphics[scale=.29]{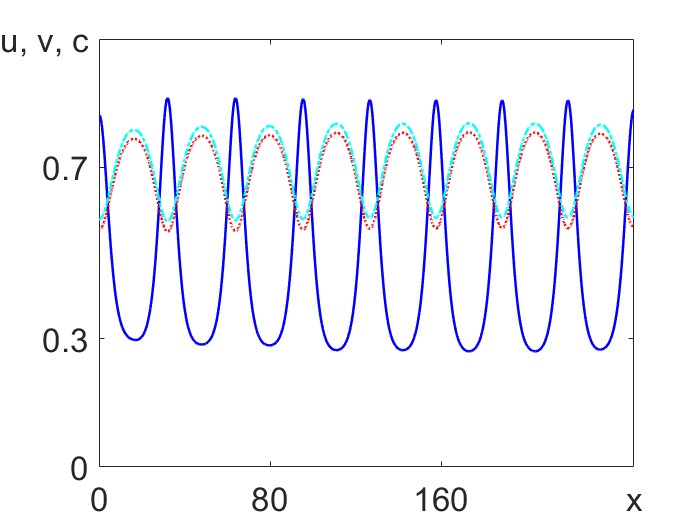}
\end{subfigure}
\caption{\em{\textbf{Stability of the coexistence steady state depending on interspecific competition}. \textbf{(a)}: Conditions (\ref{RHcriteria}) are verified for fixed parameters $D_1=D_2=1$, $r_1=r_2=0.1$, $\chi=-10$ and $k=0.2$, in the region $b_1, \mbox{ } b_2 <1$. In domain $\delta_2$, condition for $a_3$ is violated, i.e. $a_3<0$. \textbf{(b)}: stationary periodic pattern obtained from numerical simulations by fixing $b_1=b_2=0.7 \in \delta_2$. Solid blue and dotted red lines represent the density of  $u$ and $v$, respectively. Dash-dotted cyan line shows the concentration profile of the chemical which almost coincides with that of $v$.}}
\label{domains_deltas}
\end{figure}

The aim of this section is to obtain conditions for Turing-type instability of the coexistence steady state,
\begin{equation}\label{coex}
\displaystyle (u_2, \mbox{ }v_2,\mbox{ } c_2)=\left(\frac{b_1-1}{b_1b_2-1}, \mbox{ }\frac{b_2-1}{b_1b_2-1},\mbox{ } \frac{b_2-1}{b_1b_2-1}\right).
\end{equation}
It is known that in the absence of diffusion and advection this state is stable if $b_1<1$ and $b_2<1$ (weak competition). The natural step is to evaluate the characteristic matrix \eqref{Mrda} at the coexistence state and investigate under what conditions it will have positive eigenvalues. As the matrix \eqref{Mrda} is 3x3 matrix its eigenvalues are roots of the cubic polynomial: 
\begin{equation}\label{carpol}
\lambda^3+a_1\lambda^2+a_2\lambda +a_3=0,
\end{equation}
where coefficients $a_1$, $a_2$ and $a_3$ are defined by the entries of matrix \eqref{Mrda} where $(u^*,v^*,c^*)=(u_2,v_2,c_2)$. Analytical expression for these coefficients are awkwardly long to be presented here, however their numerical values can always be found for any given set of parameter values of model \eqref{twospecieschem}. The Routh-Hurwitz criteria states that if
\begin{equation}\label{RHcriteria}
\begin{cases}
a_1,\mbox{ } a_2, \mbox{ } a_3 >0,\\
a_3-a_1a_2<0,
\end{cases}
\end{equation}
then all roots of \eqref{carpol} are negative and the steady state is stable. If any of these conditions is violated, then at least one eigenvalue is positive and the steady state is unstable. For set of parameter values $D_1=D_2=1$, $r_1=r_2=0.1$, $\chi=-10$ and $k=0.2$ we have found (using MATLAB code) how the coefficients $a_1, \mbox{ } a_2, \mbox{ } a_3$ and the expression $a_3-a_1a_2$ depend on the model parameters $b_1$ and $b_2$. Fig.\ref{domains_deltas}$(a)$ shows that the region of interest, $b_1<1$, $b_2<1$, is split into two domains $\delta_1$ and $\delta_2$. In domain $\delta_1$ all conditions of the Routh-Hurwitz criteria for stability \eqref{RHcriteria} are met, so the coexistence state is stable if $(b_1, b_2) \in \delta_1$. However, in domain $\delta_2$ the Routh-Hurwitz criteria for stability are violated, namely, $a_3<0$. This means that if  $(b_1, b_2) \in \delta_2$  then the coexistence steady state is unstable and stationary periodic patterns are expected to form. This conclusion is confirmed by numerical simulations shown in Fig.\ref{domains_deltas}$(b)$. 

\begin{figure}[h]
\centering
\begin{subfigure}{.32\textwidth}
\caption{}
\includegraphics[scale=.3]{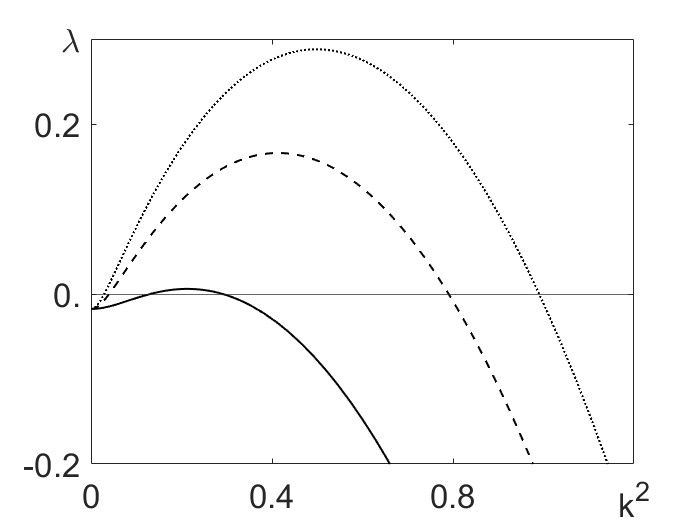}
\end{subfigure}
\hfill
\begin{subfigure}{.47\textwidth}
\caption{}
\includegraphics[scale=.3]{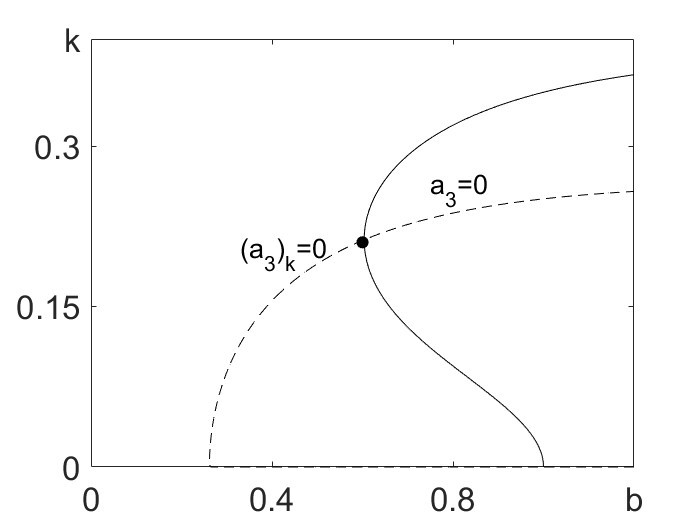}
\end{subfigure}
\caption{\em{\textbf{Most unstable wavenumber $k$ for different chemotactic strengths and competition factors.} \textbf{(a)}: The wavelength $\lambda$ as function of the wavenumber $k$ for different chemotactic strengths: $\chi=\{-10 \mbox{ (solid) }, -50 \mbox{ (dashed) }, -90 \mbox{ (dotted) }\}.$  \textbf{(b)}: Minimum value of $b$ that initiates formation of patterns and the most unstable wavenumber $k$ that ensures maximal size of $\delta_2$. Fixed parameters: $D_1=D_2=1$, $r_1=r_2=0.1$, $\chi=-10$ and $b_1=b_2=b$.  }}
\label{wavelength_wavenumber}
\end{figure}

The shape and size of the domain $\delta_2$ depend not only on the model parameters but also on the wavenumber $k$. An important question is: what is the domain $\delta_2$, encompassing all possible $(b_1, b_2)$ values that will trigger a breakdown of stability and the formation of patterns for a given set of model parameters $D_1, D_2, r_1, r_2$ and $\chi$ and any value of the wavenumber $k$? Formally, this domain is a union of domains for all possible values of $k$. However it can be estimated as a domain for the value of $k$, which gives the highest positive value for eigenvalue of matrix \eqref{Mrda}. Fig.\ref{wavelength_wavenumber}$(a)$ shows how the value of this eigenvalue depends on the wavenumber for three values of chemotactic sensitivity $\chi$. Wavenumber $k$ corresponding to the most unstable mode for each of these plots is that giving their maxima. The domain $\delta_2$ for this value of $k$ could be considered as a domain of instability for given set of parameters. As we can see from Fig.\ref{wavelength_wavenumber}$(a)$ as chemotaxtic sensitivity $\chi$ increases, the value of most unstable wavenumber also increases, resulting in smaller wavelengths and, consequently, more spikes. Additionally, as shown in the previous chapter, $\displaystyle k = \frac{i\pi}{L}$ satisfies the boundary conditions. The most unstable wavenumber can provide a rough estimate of how many spikes are expected in the system. For example, when $\chi = -10$, the most unstable wavenumber is $k = 0.2$. The number of half spikes can be predicted using $\displaystyle i = \frac{kL}{\pi}$; for a medium of length $L = 250$, $i \approx 16$, resulting in 8 full spikes. This prediction aligns well with the results shown in Fig.\ref{domains_deltas}$(b)$.

Fig.\ref{wavelength_wavenumber}$(b)$ demonstrates another method for determining the domain $\delta_2$ where the coexistence state of the system \eqref{twospecieschem} becomes unstable.  The presented analysis is performed using an assumption that $b_1 = b_2 = b$. As the loss of stability is associated with the sign change of $a_3$ (in \eqref{carpol}) we plot the line $a_3=0$ on ($b, k$) plane and note that $a_3$ is positive on the side of this plot corresponding to larger values of $b$. Minimal value of $b$ on this plot gives the the minimal value of $b$ in the domain $\delta_2$. The way to find this value is to draw another line where the partial derivative $\displaystyle \frac{\partial a_3}{\partial k} = 0$. By examining the intersection between these two curves (see Fig.\ref{wavelength_wavenumber}$(b)$), both the most unstable wavelength and the minimum value of $b$ required for pattern formation can be found. For instance, in Fig.\ref{wavelength_wavenumber}$(b)$, it can be observed that for $D_1 = D_2 = 1$, $r_1 = r_2 = 0.1$, and $\chi = -10$, the maximum size of $\delta_2$ is reached when $k \approx 0.2$, which is consistent with the results in panel $(a)$ of the same figure. Additionally, panel $(b)$ predicts that the minimum value of the competition rate $b$ needed to initiate periodic patterns is $b = 0.6$.

So far, using the Routh-Hurwitz criteria, we have identified a domain, $\delta_2$ on the plane $(b_1, b_2)$, where stationary periodic patterns are expected to emerge. Our next task is to examine how the size of this domain is affected by the remaining model parameters. We will perform this study using the technique illustrated in Fig.2$(b)$, that is by assuming that both species exert the same competitive effect on one another, i.e. $b_1 = b_2 = b$, and examining the intersection point between $\displaystyle \frac{\partial a_3}{\partial k} = 0$ and $a_3 = 0$. In what follows we will use the following default set of parameters: $D_1 = D_2 = 1$, $r_1 = r_2 = 0.1$ and $\chi = -10$.

\begin{figure}[H]
\centering
\begin{subfigure}{.16\textwidth}
\caption{}
\includegraphics[scale=.16]{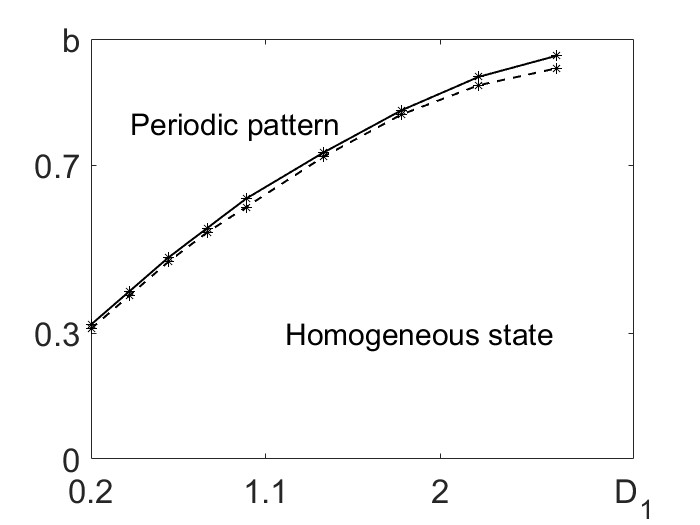}
\end{subfigure}
\hfill
\begin{subfigure}{.25\textwidth}
\caption{}
\includegraphics[scale=.16]{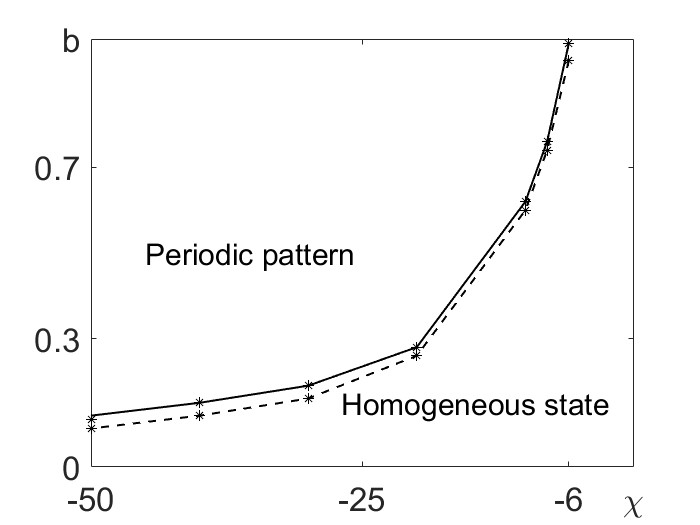}
\end{subfigure}
\begin{subfigure}{.16\textwidth}
\caption{}
\includegraphics[scale=.16]{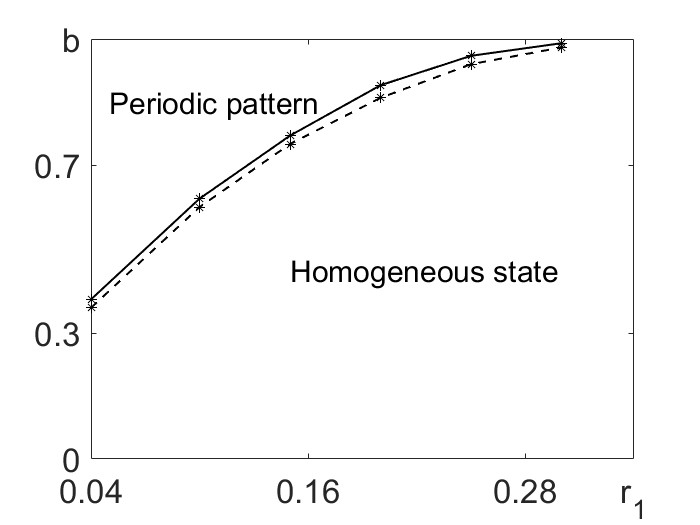}
\end{subfigure}
\hfill
\begin{subfigure}{.23\textwidth}
\caption{}
\includegraphics[scale=.16]{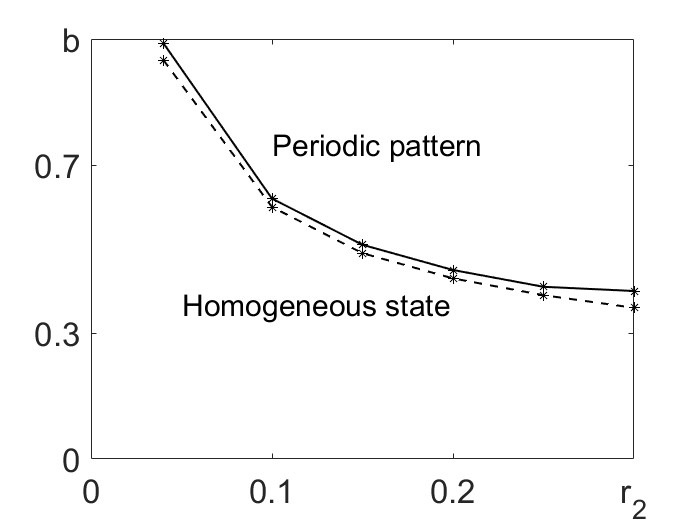}
\end{subfigure}
\caption{\em{\textbf{Domains for types of patterns that can form after disturbance of the coexistence state in the weak competition case in model \eqref{twospecieschem}. } Domains for the formation of stationary periodic patterns on the planes \textbf{(a)}: ($b,D_1$); \textbf{(b)}: ($b,\chi$); \textbf{(c)}: ($b,r_1$); \textbf{(d)}: ($b,r_2$).  Used default set of parameters: $D_1 = D_2 = 1$, $r_1 = r_2 = 0.1$ and $\chi = -10$. Solid lines represent results from numerical simulations and dashed lines represent analytical results. }}
\label{TP_b}
\end{figure}

Model (\ref{twospecieschem}) accounts for two different diffusion coefficients, $D_1$ and $D_2$, corresponding to the two species $u$ and $v$, respectively. Both analytical and numerical results show that fixing $D_1$ and varying $D_2$ has the same effect on the stability of the coexistence steady state as fixing $D_2$ and varying $D_1$. Fig.\ref{TP_b}$(a)$ shows how the minimal value of $b$ in the domain $\delta_2$ depends on the diffusion rate $D_1$ keeping the default values for other parameters. It shows a good match between the analytical and numerical results with both indicating that as diffusion increases the size of the domain $\delta_2$ decreases and it disappears when $D_1=2.6$ (shown by the vertical dashed line) and indicating that no periodic patterns will form for diffusion values greater than this. If $D > 2.6$, then only the formation of travelling waves is possible in  model (\ref{twospecieschem}). 

In a similar manner, we have investigated the effects of chemotaxis on the size of the domain $\delta_2$. Plots shown in Fig.\ref{TP_b}$(b)$
indicate that the minimal value of $b$ increases (and, therefore, the size of $\delta_2$ decreases) when the chemotactic response gets weaker ($\chi$ increases). This means that it is more likely to obtain periodic patterns for a stronger chemotactic response, while weaker chemotaxis can only result in the formation of travelling waves. Periodic patterns cannot be obtained if $\chi$ is below the threshold $\chi=-6$. Again, the shown plots indicate a good match between analytical and numerical results. 

In the last two panels of Fig.\ref{TP_b} the effects of species reproduction, $r_1$ and $r_2$, on the minimal value of $b$ in the domain $\delta_2$ are shown. We see that the reproduction rates $r_1$ and $r_2$ have opposite effects on the size of $\delta_2$. Increasing $r_1$ results to the shrinkage of $\delta_2$ with the threshold value of $r_1=0.3$ above which no periodic patterns can be obtained. On the other hand, increasing $r_2$ shows that the value of $b$ decreasing and therefore the size of the domain $\delta_2$ increasing. In this case, there is a threshold value $r_2=0.04$ below which no periodic patterns (only travelling waves) can be observed. 

Thus far, classical Turing-type instability analysis has been employed to demonstrate the formation of stationary periodic patterns caused by disturbances to the coexistence steady state of model (\ref{twospecieschem}). The Routh-Hurwitz criteria have been used to show that, for fixed parameters, a domain of instability $\delta_2$ exists, and any competition values within this domain will result in a breakdown of stability. The maximal size of this domain has been determined by identifying the most unstable wavenumber $k$, ensuring that all possible $(b_1, \mbox{ }b_2)$ values are considered, along with the minimum value of $b = b_1 = b_2$ required for pattern formation.

\subsection{Analysis of the extinction steady states}

Now we investigate whether the pattern formation can take place in system \eqref{twospecieschem} when the competion is not weak, that is when either $b_1>1$ and/or $b_2>1$. Linear analysis would point to this possibility if the stability of any of the extinction steady states, $(1, 0, 0)$ and $(0, 1, 1)$, is affected by the diffusion or advection terms. This implies that when $b_2 > 1$ periodic patterns can form if the steady state $(1, 0, 0)$ is unstable. Similarly, when $b_1 > 1$ periodic patterns can form if the steady state $(0, 1, 1)$ is driven unstable in the presence of diffusion and/or advection terms. Linear stability analysis involves examination of the eigenvalues of the characteristic matrix (\ref{Mrda}) evaluated at each extinction steady state to  identify conditions under which at least one eigenvalue has a positive real part. The analysis of the steady state $(1, 0, 0)$ is done by investigating the eigenvalues of the characteristic matrix (\ref{Mrda}) evaluated at $(n, v, c)=(1, 0, 0)$, i.e.:
\begin{equation}\displaystyle M_{(1, 0, 0)}=
\begin{pmatrix}
\displaystyle -D_1k^2-r_1 & -r_1b_1 & -\chi k^2 \\
\displaystyle 0 & -D_2k^2-r_2(b_2-1) & 0 \\
\displaystyle 0 & 1& -1-k^2
\end{pmatrix},
\end{equation}
which has 3 eigenvalues:
\begin{equation}
\begin{cases}
\displaystyle -1-k^2,\\
\displaystyle -D_1k^2-r_1,\\
\displaystyle -D_2k^2-r_2(b_2-1),
\end{cases}
\end{equation}
with negative real parts if $b_2>1$. This means that the steady state $(1, 0, 0)$ remains stable when perturbed in the full reaction-diffusion-advection system and we should not expect formation of stationary patterns from this state. This result also holds for the other extinction steady state,  $(0, 1, 1)$: if $b_1>1$ this steady state is stable irrespective on the strengths of both diffusions and chemotaxis.  

\begin{figure}[h]
	\centering
	\begin{subfigure}{.16\textwidth}
		\caption{}
		\includegraphics[scale=.16]{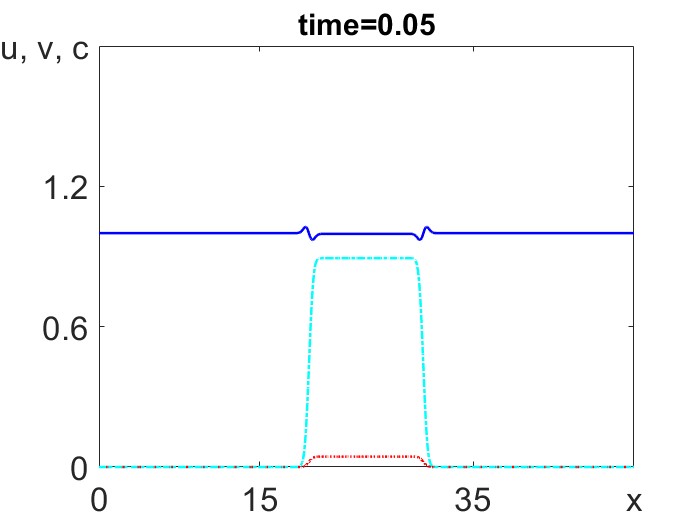}
	\end{subfigure}
	\hfill
	\begin{subfigure}{.25\textwidth}
		\caption{}
		\includegraphics[scale=.16]{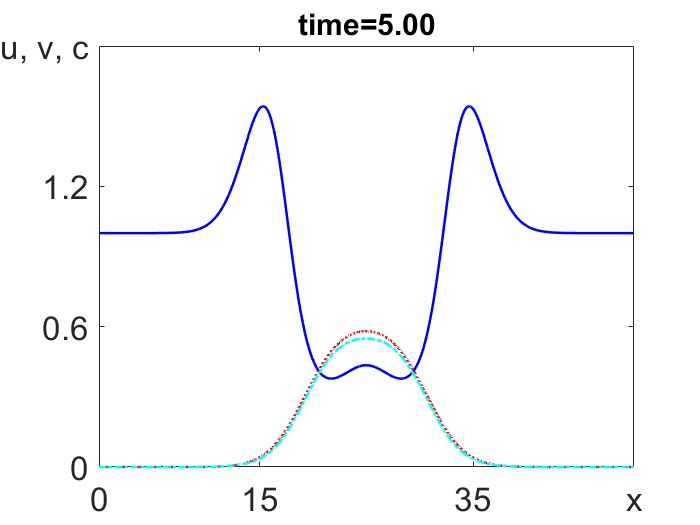}
	\end{subfigure}
	\begin{subfigure}{.16\textwidth}
		\caption{}
		\includegraphics[scale=.16]{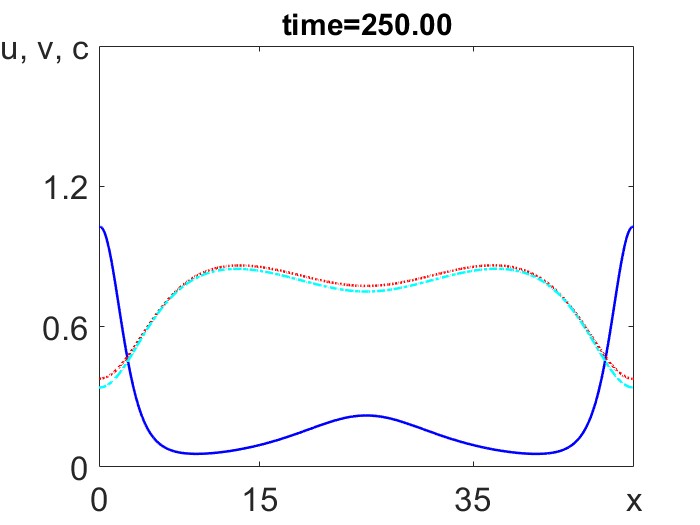}
	\end{subfigure}
	\hfill
	\begin{subfigure}{.23\textwidth}
		\caption{}
		\includegraphics[scale=.16]{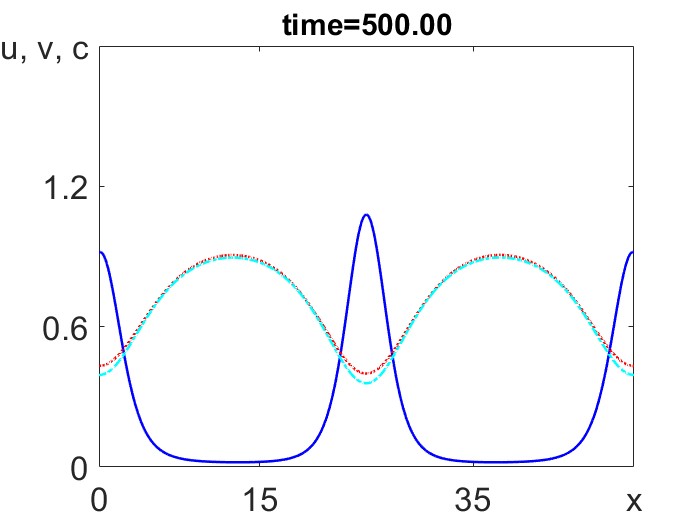}
	\end{subfigure}
	\caption{\em{\textbf{Simulation of pattern obtained in the weak-strong competition case when $b_2>1$.} \textbf{(a)}: An initially large perturbation of the steady state $(1, 0, 0)$ (dash-dotted cyan line). \textbf{(b)}: $v$ starts producing a chemical agent $c$ (dotted red line) which repels species $u$ (solid blue line). \textbf{(c)}: Density of $u$ slowly starts increasing in the middle, while the density of $v$ and concentration of $c$ slowly start decaying \textbf{(d)}: Stationary pattern consisting of two full spikes has formed. Parameter values: $D_1=D_2=1, \mbox{ } r_1=r_2=0.1, \mbox{ } \chi=-10$, $b_1=0.7, \mbox{ } b_2=1.7$, medium size, $L=50$, and the amplitude of perturbation, $\tilde{v}=0.9.$ }}
	\label{twosp_100_sim}
\end{figure}

One characteristic of the periodic patterns obtained so far is that they have formed from infinitesimal perturbations, meaning that as long as the conditions for the breakdown of stability have been met, the size of the perturbation did not matter and has not been a factor to consider. However, when investigating the steady state $(1, 0, 0)$, one question to consider is whether a large enough perturbation in the second species, $v$, would result in the formation of patterns, or if the perturbation would decay and the system would return to being homogeneous. The biological idea behind this is that in experiments, if the starting point were a Petri dish with a species of bacteria, $u$, distributed homogeneously, introducing a large enough quantity of the second species of bacteria, $v$, could produce enough chemotactic agent to start repelling the first species, and with the aid of competition, the second species could survive and reproduce. The repulsion of the first species, $u$, by the chemical, $c$, could launch the aggregation and hence, the formation of spots. In order to test this hypothesis, computational simulations have been performed with fixed parameters: $D_1=D_2=1$, $\chi=-10$, $r_1=r_2=0.1$, $b_1=0.7$, $b_2=1.7$, such that $b_2>1$ and $v=0+\tilde{v}$, where the perturbation $\tilde{v}=0.9$. Simulated profiles are shown in Figure \ref{twosp_100_sim}.

Figure \ref{twosp_100_sim} shows that, contrary to linear analysis, stationary periodic patterns can emerge from finite perturbation of the stable extinction steady state $(1, 0, 0)$. An important characteristic of this pattern is that aggregation can only be initiated by the perturbation of species $v$. Perturbing species $u$, or the chemical $c$, results in the system decaying back to homogeneity, regardless of the amplitude of the perturbation. This leads us to the next question, which is how the minimal amplitude required for pattern formation is affected by changes in the model parameters.

Next, with the aid of computational simulations, the effect of model parameters on the minimal perturbation amplitude required for pattern formation is investigated. The default set of model parameters used in these simulations are: $D_1=D_2=1$, $r_1=r_2=0.1$, $\chi=-10$, $b_1=0.7$, $b_2=1.7$ and medium size, $L=50$. By varying these parameters one at a time, we note the minimal amplitude required for a breakdown of stability. It is also important to note that, unlike patterns around the coexistence steady state, travelling waves do not emerge from the extinction steady state. This means that, depending on the amplitude of the perturbation, the system can exhibit pattern formation or return to a homogeneous state. 

\begin{figure}[H]
	\centering
	\begin{subfigure}{.3\textwidth}
		\caption{}
		\includegraphics[scale=.2]{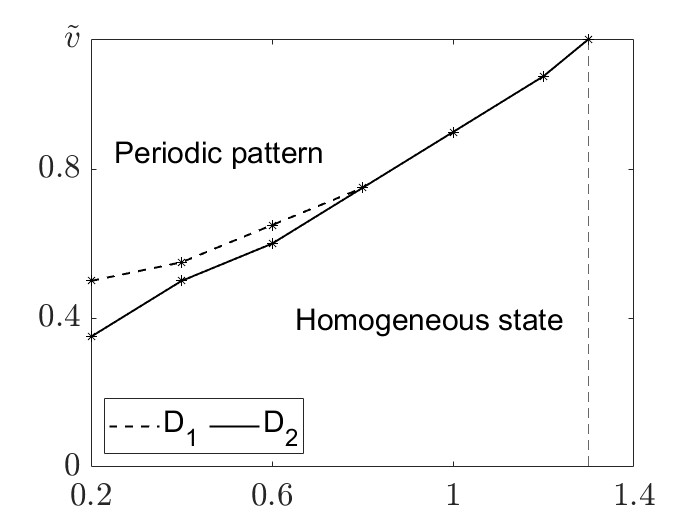}
	\end{subfigure}
	\hfill
	\begin{subfigure}{.33\textwidth}
		\caption{}
		\includegraphics[scale=.2]{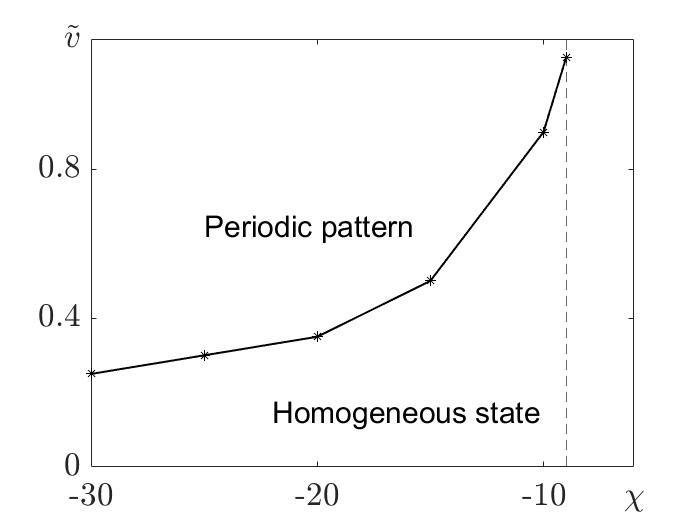}
	\end{subfigure}
	\begin{subfigure}{.3\textwidth}
		\caption{}
		\includegraphics[scale=.2]{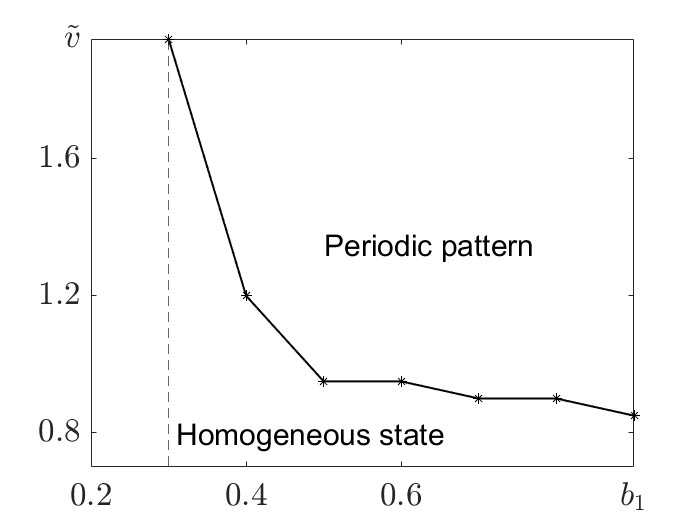}
	\end{subfigure}
	\hfill
	\begin{subfigure}{.3\textwidth}
		\caption{}
		\includegraphics[scale=.2]{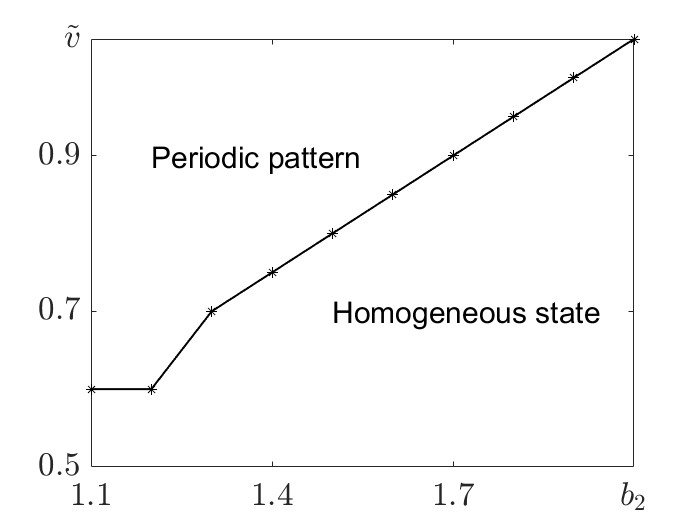}
	\end{subfigure}
	\hfill
	\begin{subfigure}{.3\textwidth}
		\caption{}
		\includegraphics[scale=.2]{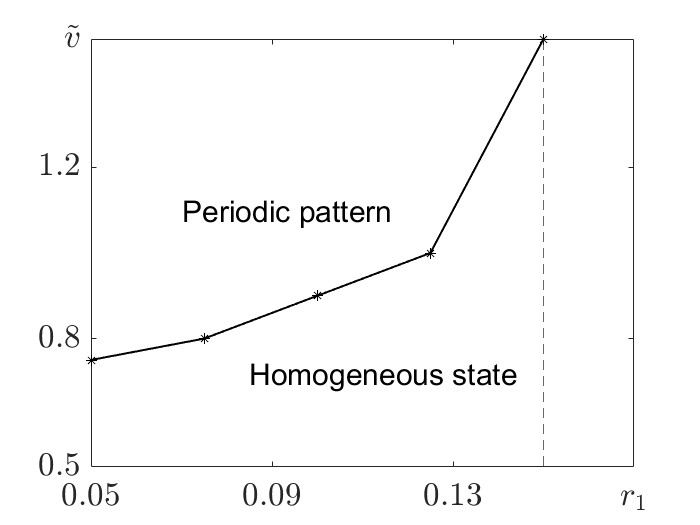}
	\end{subfigure}
	\hfill
	\begin{subfigure}{.3\textwidth}
		\caption{}
		\includegraphics[scale=.2]{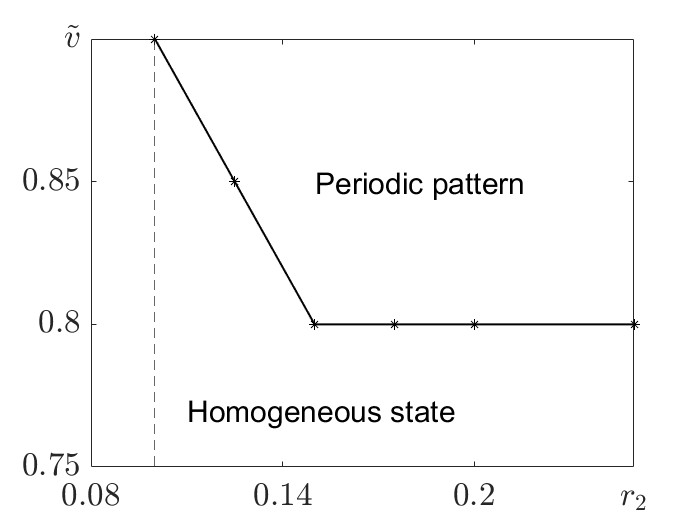}
	\end{subfigure}
	\caption{\em{\textbf{The dependence of the minimal perturbation amplitude, $\tilde{v}$, from the steady state $(1, 0, 0)$, required to generate periodic patterns, on the model parameters.} Dependence on the diffusion coefficients (panel $(a)$),  chemotactic sensitivity (panel $(b)$), competition strength $b_1$ (panel $(c)$), and $b_2$ (panel $(d)$), reproduction rate $r_1$ (panel $(e)$) and $r_2$ (panel $(f)$). Default set of parameters: $D_1=D_2=1$,  $\chi=-10$, $r_1=r_2=0.1$, $b_1=0.7$ and $b_2=1.7$. }}
	\label{pert_param}
\end{figure}

The plots in Fig.\ref{pert_param}$(a)$ show the minimum perturbation amplitude required in simulations for periodic patterns to emerge from the extinction steady state $(1, 0, 0)$ when both diffusion coefficients are varied one at a time and the other parameters are fixed. We see that the impacts of diffusion coefficients are similar  with the increase of the threshold perturbation, $\tilde{v}$, with the increase of the diffusion. 
An important observation is that for $D_1 > 1.3$ or $D_2 > 1.3$, the system does not show the formation of patterns, meaning that from any perturbation the system relaxes back to homogeneous state. 

The effect of chemotaxis on the perturbation amplitude is shown by the plot in Fig.\ref{pert_param}$(b)$. For chemoattraction, or weak chemorepulsion, when $\chi>-9$, the system cannot be perturbed from the homogeneous state, no matter how large the perturbation is. However, for stronger chemorepulsion, when $\chi < -9$, the system can be transformed into inhomogeneous state provided the disturbance is large enough. As chemorepulsion increases, a smaller perturbation amplitude is needed to generate pattern formation. 

Similarly, the effects of reproduction rates on the threshold perturbation amplitude are also investigated with the results presented in Fig.\ref{pert_param}$(c,d)$. It appears that the two reproduction coefficients, $r_1$ and $r_2$, corresponding to species $u$ and $v$, respectively, have opposite effects on the threshold amplitude. If the reproduction rate of the species $u$ is high ($r_1 > 0.15$), then the system cannot be transformed into stationary inhomogeneous state no matter how large is perturbation. Additionally, as $r_1$ decreases, the minimal perturbation amplitude required for the formation of stationary periodic pattern decreases, meaning that if $u$ reproduces more slowly, it is more likely to initiate aggregation through perturbation. On the other hand, if $r_2 < 0.1$, no pattern formation in the system can be observed. Moreover, as $r_2$ increases, the amplitude required for a breakdown of stability decreases, converging to $\tilde{v} = 0.8$ for $r_2 \geq 0.15$. 

The effect of the competition rates, $b_1$ and $b_2$, on the threshold perturbation amplitude is presented by plots in Fig.\ref{pert_param}$(e,f)$. Numerical results show that pattern can form only if $b_1$ is large enough, $b_1 \geq 0.3$. if $b_1 < 0.3$, the system relaxes back to a homogeneous state from any initial state. Conversely, as $b_2$ increases, the perturbation amplitude required for instability also increases. This means that, in order to initiate pattern formation, higher competition of $u$ on $v$, represented by $b_2$, requires a larger disturbance for species $v$ to survive and reproduce.

In this section, we have shown that, contrary to linear analysis, which suggests that the steady state $(1, 0, 0)$ of system (\ref{twospecieschem}) remains stable if $b_2 > 1$, stationary inhomogeneous pattern can emerge in the system \eqref{twospecieschem} from a finite amplitude disturbance of this steady state. This phenomenon has been demonstrated in numerical simulations. Using numerical simulations, the effect of model parameters on the minimum perturbation amplitude required to generate a breakdown of stability and, consequently, pattern formation has been investigated. We have also found that if the strong-weak competition case when $b_2<1$ and $b_1>1$ we don't observe pattern formation no matter what initial conditions are set.

While linear analysis focuses on the breakdown of stability and formation of inhomogeneous patterns, Fourier analysis can be used to describe the pattern once it has formed and stabilized. In particular, Fourier analysis can be applied to obtain information about the pattern's characteristics, such as amplitude and wavelength.

\section{Fourier analysis of patterns obtained in numerical simulations}

In this section we will focus on Fourier series representing stationary periodic patterns obtained in numerical simulations, like those shown in Fig.\ref{domains_deltas}$(b)$ and Fig.\ref{twosp_100_sim}$(d)$. As the stationary pattern is that which does not change over time, it satisfies the simplified version of model (\ref{twospecieschem}):
\begin{equation}\displaystyle\label{twospeciestationary}
\displaystyle\begin{cases}
\displaystyle D_1n_{xx}-\chi(uc_x)_x+r_1u(1-u-b_1v)=0,\\
\displaystyle D_2v_{xx}+r_2v(1-v-b_2u)=0,\\
\displaystyle c_{xx}+v-c=0,
\end{cases}
\end{equation}
and can be represented by Fourier series:
\begin{equation} \label{twospecies_Fouriersol}
\begin{cases}
\displaystyle u=\sum_{i=0}^M \alpha_i\cos\left(\frac{i\pi}{L}x\right),\\[7pt]
\displaystyle v=\sum_{i=0}^M \gamma_i\cos\left(\frac{i\pi}{L}x\right),\\[7pt]
\displaystyle c=\sum_{i=0}^M \beta_i\cos\left(\frac{i\pi}{L}x\right).
\end{cases}
\end{equation}
The coefficients $\alpha_i$, $\gamma_i$ and $\beta_i$ define the amplitudes of mode $i$ for the variables $u$, $v$ and $c$. For smooth profiles, these coefficients quickly tend to zero as $i$ increases, allowing us to truncate the series by considering only the first $M$ terms, where $M$ should be carefully chosen. For a known profile, $u(x)$, the coefficients are found using the formulas:
\begin{equation*}
\alpha_0=\frac{1}{L}\int_{0}^{L}u(x)dx, 
\end{equation*}
and for $i>0$
\begin{equation} \label{alphas}
\alpha_i=\frac{2}{L}\int_{0}^{L}u(x)\cos\frac{i\pi x}{L}dx.  
\end{equation}
In the rest of this chapter, we will focus on profiles $u(x)$ (and coefficients $\alpha_i$), keeping in mind that the analysis of the profiles $v(x)$ and $c(x)$ is done in the same way.

\begin{figure}[h]
\centering
\begin{subfigure}{.25\textwidth}
\caption{}
\includegraphics[scale=.21]{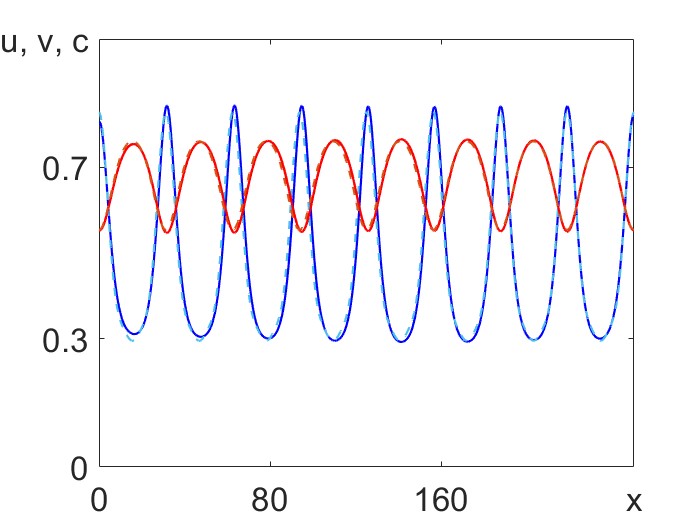}
\end{subfigure}
\hfill
\begin{subfigure}{.25\textwidth}
\caption{}
\includegraphics[scale=.21]{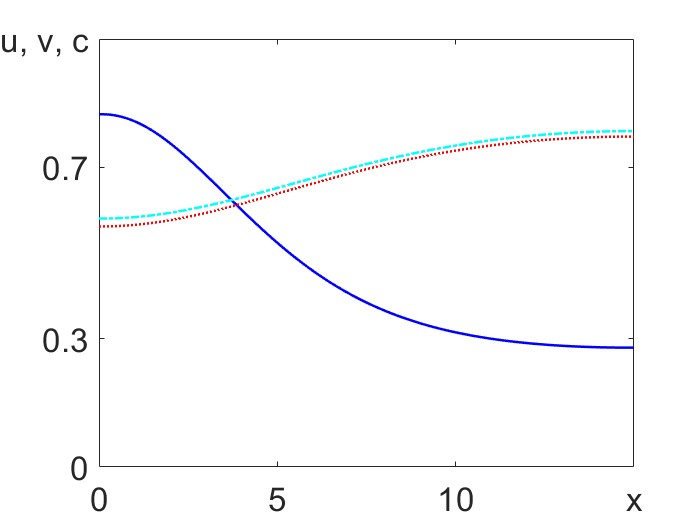}
\end{subfigure}
\hfill
\begin{subfigure}{.32\textwidth}
\caption{}
\includegraphics[scale=.21]{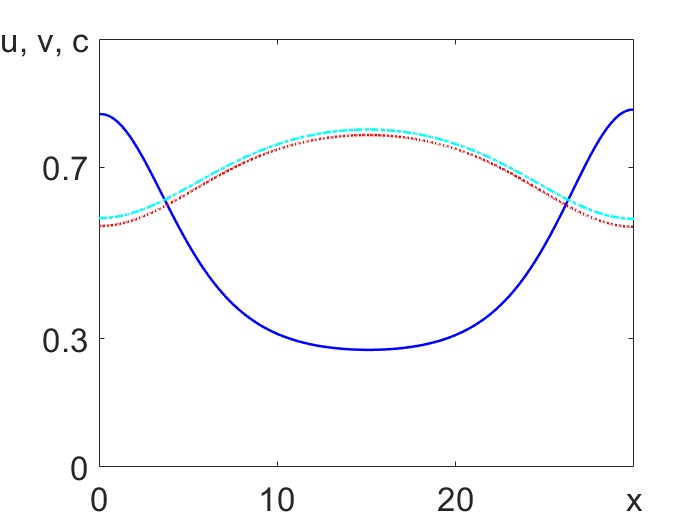}
\end{subfigure}
\caption{\em{\textbf{Numerical simulations of the stationary patterns formed in model (\ref{twospecieschem}) for different medium lengths. } \textbf{(a)}: Comparison between the numerical profile from simularions (solid) and the profile including only the highest three modes (dashed) given in \eqref{coeffs16} for a large medium  $L=250$. \textbf{(b)}: Reduced medium length such that only half a spike produced $L=15$. \textbf{(c)}: Full spike for medium legth $\L=30$. Blue and red lines represent the densities of $u$ and $v$, respectively, while cyan lines the concentration of the chemical $c$. Used  parameter values: $D_1=D_2=1$, $r_1=r_2=0.1$, $\chi=-10$ and $b_1=b_2=0.7$.}}
\label{simulations_twosp}
\end{figure}

Formulas \eqref{alphas} can be used for the spectral decomposition of patterns obtained numerically. A typical stationary solution obtained from numerical simulations of the system \eqref{twospecieschem} is shown in Fig.\ref{domains_deltas}$(b)$. Spectral decomposition of this profile reveals that only four modes have reasonably high coefficients:


\begin{equation}  \label{coeffs16}
\begin{cases}
\displaystyle \alpha_0=0.4777, \mbox{ } \alpha_{16}=0.2445, \mbox{  } \mbox{ } \alpha_{32}=0.0842, \mbox{ } \mbox{ } \mbox{ } \alpha_{48}=0.0225;\\
\displaystyle \gamma_0=0.6719, \mbox{ } \gamma_{16}=-0.1024, \mbox{ } \gamma_{32}=-0.0145, \mbox{ } \gamma_{48}=-0.0010;\\
\displaystyle \beta_0=0.6719, \mbox{ } \beta_{16}=-0.0984, \mbox{ } \beta_{32}=-0.0124, \mbox{ } \beta_{48}=-0.0006,  
\end{cases}
\end{equation}
while all other coefficients are considerably smaller (all other $\alpha$-s are less than 0.01). $\alpha_0$ represents the average level of $u$ for the entire pattern, $\alpha_{16}$ defines the amplitude of the mode with a characteristic length of 1/16th of the domain size, which corresponds to the 8 spikes observed. This matches the number of spikes seen in Fig.\ref{domains_deltas}$(b)$. Finally, $\alpha_{32}$ and $\alpha_{48}$ correspond to the second and third harmonics of the main harmonic given by $\alpha_{16}$. Thus, the amplitude of the pattern shown in Fig.\ref{domains_deltas}$(b)$ can be estimated as $2\alpha_{16}$. Varying the size of the domain will change the number of observed spikes but will not affect their amplitude or spatial periodicity. This is illustrated by the patterns shown in panels $(b)$ and $(c)$ of Fig.\ref{simulations_twosp}, which are obtained in simulations of smaller domains where only half $(b)$ and one $(b)$ spike can fit.

\begin{figure}[h]
\centering
\begin{subfigure}{.31\textwidth}
	\caption{}
	\includegraphics[scale=.31]{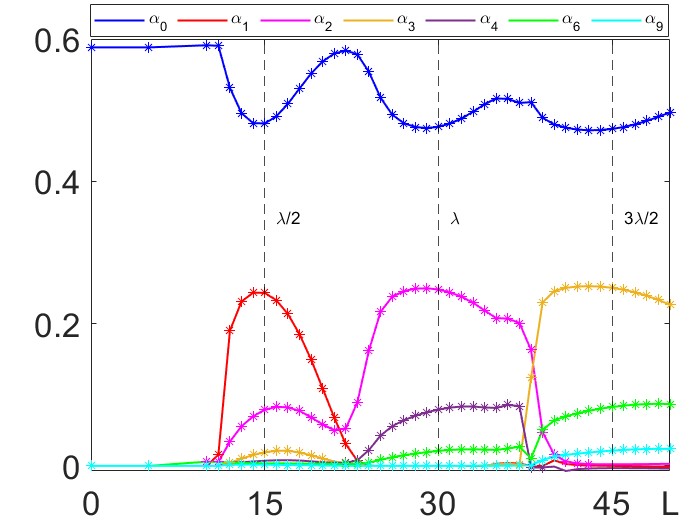}
\end{subfigure}
\hfill
\begin{subfigure}{.48\textwidth}
	\caption{}
	\includegraphics[scale=.31]{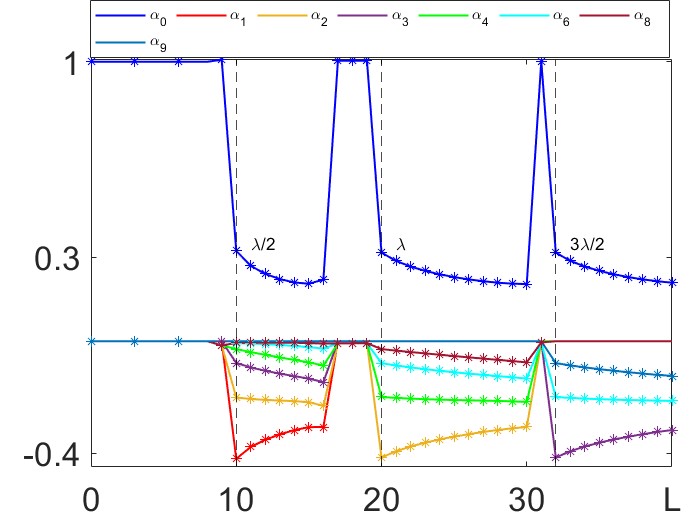}
\end{subfigure}
\caption{\em{\textbf{Dependence of Fourier coefficients $\alpha_i$ on medium size $L$ for numerically simulated  $u$-profiles in model \eqref{twospecieschem}.} For patterns forming \textbf{(a)}: in weak competition case, $b_2=0.7$ and \textbf{(b)} in weak-strong competition case, $b_2=1.7$.   Other model parameter: $D_1=D_2=1$, $r_1=r_2=0.1$, $\chi=-10$,  $b_1=0.7$.}}
\label{ag_numerically}
\end{figure}

Patterns occurring in a small domain (like those shown in Fig.\ref{simulations_twosp}$(b)$ and $(c)$) are of particular interest, as the corresponding Fourier series can be truncated at a reasonably low value of $M$ (see Eq. \eqref{twospecies_Fouriersol}). Using Fourier decomposition of numerically obtained patterns, we investigated how the coefficients of Fourier modes depend on the domain size. In Fig. \ref{ag_numerically}, we show plots of $\alpha_i$ for $i=0, 1, 2, 3, 4, 6, 9$ against the domain size, $L$, which varies from $0$ to $50$. We observe that while $\alpha_0$ does not vary much (staying in the range between $0.5$ and $0.6$) and the omitted $\alpha$-coefficients are always negligibly small, the displayed coefficients vary significantly. For $L < 10$, $\alpha_0 = 0.6$ and all other coefficients are zero, indicating that the system is in a homogeneous state with no spikes formed. For $10 < L < 20$, the first coefficient, $\alpha_1$, is larger than any subsequent coefficient ($\alpha_1 > \alpha_i$, $\forall i > 1$), reflecting the fact that only half of a spike can fit in. Within this range of domain sizes, the value of $\alpha_1$ increases from zero, reaches a maximum value $\alpha_1 (= \alpha_{\text{max}}) = 0.2445$ at $L = 15$, and then decreases to zero. This domain size, which corresponds to the maximal value of $\alpha_1$, will be considered the characteristic length (or half-wavelength) of the periodic pattern and will be denoted as $\Lambda_0$, i.e., $\Lambda_0 = 15$ for the system \eqref{twospecieschem} with the model parameters used in simulations to produce Fig. \ref{ag_numerically}.

As we seen earlier, stationary periodic patterns in the system \eqref{twospecieschem} form under two regimes: in the case of weak competition ($b_1<1$ and $b_2<1$, see Fig.\ref{domains_deltas}$(b)$ and Fig.\ref{simulations_twosp}) and in the case of weak-strong or strong competition when $b_2>1$. In the latter case despite the stability of homogeneous state $(1,0,0)$ we observe the formation of stationary periodic patterns (see Fig.\ref{twosp_100_sim}). Characteristics of the patterns forming in these two regimes are significantly different. Fig.\ref{ag_numerically} $(b)$ shows how the first 9 $\alpha$-coefficients describing $u$-profile in the pattern obtained for $b_2>1$ depend on the medium size. As we can see the spectral decomposition of the profiles in panel $(b)$ is considerable different from those shown in the panel $(a)$. For $L \in [10, 17)$, $\alpha_1$ is the largest mode and in this case only half of the spike can fit into the simulated medium. We note that $\alpha_0$ drops significantly when pattern forms: for $L \in [17, 20)$, $\alpha_0 = 1$ is the only non-zero coefficient, meaning the system relaxes back to a homogeneous steady state. For $L \in [20, 30]$, $\alpha_2$ is the highest frequency mode, and patterns emerging for medium lengths in this region will exhibit one full spike. Additionally, the fastest-growing mode corresponds to $L = 20$, as $\alpha_2$ is largest at this point.

So far, we have shown that we can use Fourier decomposition to obtain information about the amplitude and wavelength of a stationary periodic pattern obtained in numerical simulations and express the solution of system (\ref{twospecieschem}) as a Fourier series. Our next aim is to to solve this system analytically by finding out the Fourier coefficients of solution to the system \eqref{twospeciestationary}.

\section{Analytical solution represented by Fourier series}

The analytical solution of the system (\ref{twospeciestationary}) can be found using substitution (\ref{twospecies_Fouriersol}), which transforms this system into:
\begin{equation}\label{Fseriessystem}
\begin{cases}
\displaystyle -D_1\sum_{i=0}^M (ik)^2\alpha_i\cos(ikx)+\chi\frac{\partial}{\partial x} \left(\sum_{i=0}^M\alpha_i\cos(ikx)\sum_{i=0}^M ik\beta_i\sin(ikx)\right)=0\\
\hfill +r_1\sum_{i=0}^M \alpha_i\cos(ikx)\left(1-\sum_{i=0}^M\alpha_i\cos(ikx)-b_1\sum_{i=0}^M\gamma_i\cos(ikx)\right),\\[10pt]
\displaystyle  -D_2\sum_{i=0}^M (ik)^2\gamma_i\cos(ikx)+r_2\sum_{i=0}^M \gamma_i\cos(ikx)\left(1-\sum_{i=0}^M\gamma_i\cos(ikx)-b_2\sum_{i=0}^M\alpha_i\cos(ikx)\right)=0,\\[10pt]
\displaystyle -\sum_{i=0}^M(ik)^2\beta_i\cos(ikx)+\sum_{i=0}^M\gamma_i\cos(ikx)-\sum_{i=0}^M\beta_i\cos(ikx)=0,
\end{cases}
\end{equation}
where $k=\pi/L$. 
Truncation of the system \eqref{Fseriessystem} at $M=0$ gives the system of three algebraic equations:
\begin{equation}\label{a0b0g0system}
\begin{cases}
\displaystyle r_1\alpha_0(1-\alpha_0-b_1\gamma_0)=0,\\
\displaystyle r_2\gamma_0(1-\gamma_0-b_2\alpha_0)=0,\\
\displaystyle \gamma_0-\beta_0=0,
\end{cases}
\end{equation}
which has four solutions: 
\[\displaystyle (\alpha_0, \mbox{ } \gamma_0, \mbox{ } \beta_0)=\left\{(0, 0, 0), \mbox{ } (1, 0, 0), \mbox { } (0, 1, 1), \mbox{ } \left(\frac{b_1-1}{b_1b_2-1}, \frac{b_2-1}{b_1b_2-1}, \frac{b_2-1}{b_1b_2-1}\right)  \right\},\]
corresponding to the steady states of model (\ref{twospecieschem}), as expected. Note, that the third equation in the system \eqref{Fseriessystem} gives:  
\begin{equation}\label{bigi}
\displaystyle -(ik)^2\beta_i+\gamma_i-\beta_i=0 \Longrightarrow \beta_i=\frac{\gamma_i}{1+(ik)^2},
\end{equation}
which allows to reduce the system to equations in terms of unknowns $\alpha_i$ and $\gamma_i$. Truncating the system  (\ref{Fseriessystem}) at $M > 0$ we get $2M$ simultaneous equations that need to be solved:
\begin{equation}\label{Fsimeqs_twosp}
\begin{cases}
\displaystyle \alpha_0-\alpha_0^2-\sum_{i=1}^M \frac{\alpha_i^2}{2} -b_1\left(\alpha_0\gamma_0+\sum_{i=1}^M \frac{\alpha_i\gamma_i}{2}\right)=0,\\[10pt]
\displaystyle \gamma_0-\gamma_0^2-\sum_{i=1}^M\frac{\gamma_i^2}{2} - b_2 \left(\alpha_0\gamma_0+\sum_{i=1}^M\frac{\alpha_i\gamma_i}{2}\right)=0,\\[10pt]
\displaystyle -D_1\alpha_1k^2+\chi k^2\left(\alpha_0\beta_1+\sum_{i=1}^{M-1}\frac{i+1}{2}\alpha_i\beta_{i+1}-\sum_{i=2}^M\frac{i-1}{2}\displaystyle \alpha_i\beta_{i-1}\right)-\\[10pt]
\displaystyle \hfill r_1b_1\left(\alpha_0\gamma_1+\alpha_1\gamma_0+\sum_{i=1}^{M-1}\left(\frac{\alpha_i\gamma_{i+1}}{2}+\frac{\alpha_{i+1}\gamma_i}{2}\right)\right)+\\[10pt]
\displaystyle \hfill r_1\left(\alpha_1-2\alpha_0\alpha_1-\sum_{i=1}^{M-1}\alpha_i\alpha_{i+1}\right)=0,\\[10pt]
\displaystyle -D_2\gamma_1k^2-r_2b_2\left(\alpha_0\gamma-1+\alpha_1\gamma_0+\sum_{i=1}^{M-1}\left(\frac{\alpha_i\gamma_{i+1}}{2}+\frac{\alpha_{i+1}\gamma_i}{2}\right)\right)+\\[10pt]
\displaystyle \hfill r_2\left(\gamma_1-2\gamma_0\gamma_1-\sum_{i=1}^{M-1}\gamma_i\gamma_{i+1}\right)=0,\\[10pt]
\displaystyle -4D_1\alpha_2k^2+\chi k^2 \left(4\alpha_0\beta_2+\alpha_1\beta_1+\sum_{i=1}^{M-2}(i+2)\alpha_1\beta_{i+2}-\sum_{i=3}^M(i-2)\alpha_i\beta_{i-2}\right)-\\[10pt]
\displaystyle \hfill r_1b_1\left(\alpha_0\gamma_2+\alpha_2\gamma_0+\frac{\alpha_1\gamma_1}{2}+\sum_{i=1}^{M-2}\left(\frac{\alpha_i\gamma_{i+2}}{2}+\frac{\alpha_{i+2}\gamma_i}{2}\right)\right)+\\[10pt]
\displaystyle \hfill r_2\left(\alpha_2-\frac{\alpha_1^2}{2}-2\alpha_0\alpha_2-\sum_{i=1}^{M-2}\alpha_i\alpha_{i+2}\right)=0,\\[10pt]
\displaystyle -4D_2\gamma_2k^2-r_2b_2\left(\alpha_0\gamma_2+\alpha_2\gamma_0+\frac{\alpha_1\gamma_1}{2}+\sum_{i=1}^{M-2}\left(\frac{\alpha_i\gamma_{i+2}}{2}+\frac{\alpha_{i+2}\gamma_i}{2}\right)\right) +\\[10pt]
\displaystyle \hfill r_2\left(\gamma_2-\frac{\gamma_1^2}{2}-2\gamma_0\gamma_2-\sum_{i=1}^{M-2}\gamma_i\gamma_{i+2}\right)=0,\\[10pt]
\displaystyle . . .
\end{cases}
\end{equation}
where the first two equations represent the balance for coefficients of $\cos(0)$, the third and fourth equations - for the coefficients of $\cos(kx)$, the following two equations - the coefficients of $\cos(2kx)$ and so on up to coefficients of $\cos(Mkx)$. 

Analytical solutions expressed in terms of model parameters for the simultaneous system \eqref{Fsimeqs_twosp} cannot be found for $M \geq 1$, but numerical values of $\alpha_i$, $\gamma_i$, and $\beta_i$ can be determined with the aid of Maple for any given set of model parameters $D_1$, $D_2$, $\chi$, $r_1$, $r_2$, $b_1$, and $b_2$. This method can be used to compare the Fourier series coefficients obtained analytically by solving system (\ref{Fsimeqs_twosp}) with those obtained through Fourier decomposition of the simulated profile, as presented in Figure \ref{ag_numerically}. To quantify the difference between the profiles obtained from simulations and those predicted analytically, the error between the two curves is calculated using the formula:
\begin{equation}\label{erfor}
\displaystyle ER = \int_0^L (P_A - P_N)^2dx,
\end{equation}
where $P_A$ represents the profile found analytically as a solution of the system \eqref{Fsimeqs_twosp} and $P_N$ - the profile from numerical simulations. 

Figure \ref{ag_numerically} shows that the half-wavelength of the periodic pattern obtained numerically using parameter values: $D_1=D_2=1$, $r_1=r_2=0.1$, $\chi=-10$, and $b_1=b_2=0.7$ is equal to $15$. In addition, this figure indicates that for an accurate approximation of the profile $u(x)$, coefficients up to and including $\alpha_3$ need to be considered. For the fixed parameters mentioned, we have used Maple to find coefficients up to and including $\alpha_3$ (to describe the profile $u(x)$) and $\gamma_3$ (to describe the profile $v(x)$) by solving the system (\ref{Fsimeqs_twosp}) for the medium of size, $L=15$, using truncations at $M=1, 2, 3$ . Found solutions for $\alpha$-s are presented in Table \ref{Maple_a4}.

\begin{table}[h]
\begin{center}
    \begin{tabular}{| l | l | l | l | l | l | } 
    \hline
$L=15$ & $\alpha_0$  & $\alpha_1$ & $\alpha_2$ & $\alpha_3$   \\[7pt] \hline
Numerical & 0.4813 & 0.2431 & 0.0787 & 0.0191   \\[7pt] \hline
$M=1$ & 0.1878 & 0.3330 & 0 & 0  \\[7pt] \hline
$M=2$ & 0.4753 & 0.2446 & 0.0961 & 0 \\[7pt] \hline
$M=3$ & 0.4702 & 0.2532 & 0.0849 & 0.0241  \\[7pt] \hline
    \end{tabular}
\end{center}
\caption{\em{\textbf{First four Fourier coefficients describing $u$-profile in the case of weak competition}. Numerical coefficients are obtained by Fourier decomposition of the profile $u(x)$ obtained in simulations of system \eqref{twospecieschem} and shown in Fig.\ref{simulations_twosp}(b), while analytical coefficients are obtained by solving system (\ref{Fsimeqs_twosp}) truncated at $M=1, 2$ and $3$. Model parameters: $D_1=D_2=1, \mbox{ } r_1=r_2=0.1, \mbox{ } \chi=-10, \mbox{ }b_1=b_2=0.7$ and medium length $L=15$.}}
\label{Maple_a4}
\end{table}

Analysing the coefficients presented in Table \ref{Maple_a4}, one can clearly see that truncating at $M=1$, including only $\alpha_0$ and $\alpha_1$, does not accurately reproduce the corresponding coefficients found from numerical simulations. However, increasing the truncation to $M=2$, such that the simultaneous system is extended to solve for $\alpha_0$, $\alpha_1$, and $\alpha_2$, produces a much better match to Fourier decomposition of simulated profiles. Our next aim is to determine how many $\alpha$ coefficients are needed to accurately reproduce the pattern shown by the curve $u(x)$. It is difficult to say whether increasing the truncation further, with $M>2$, would result in a better match between numerical and analytical results. To quantify the discrepancy and better visualise the difference between the numerical and analytical profiles, the numerical profile is compared against the four analytical profiles obtained from different truncations, and the error between the two curves is calculated according to (\ref{erfor}). 

\begin{figure}[h]
\centering
\begin{subfigure}{.27\textwidth}
\caption{}
\includegraphics[scale=.21]{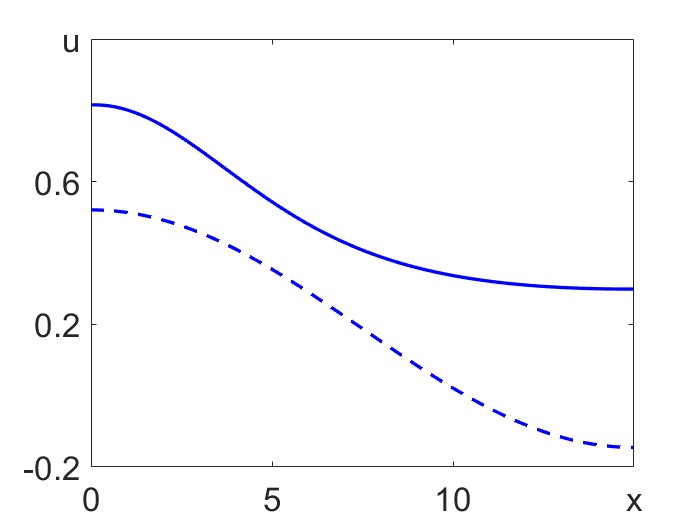}
\end{subfigure}
\hfill
\begin{subfigure}{.35\textwidth}
\caption{}
\includegraphics[scale=.21]{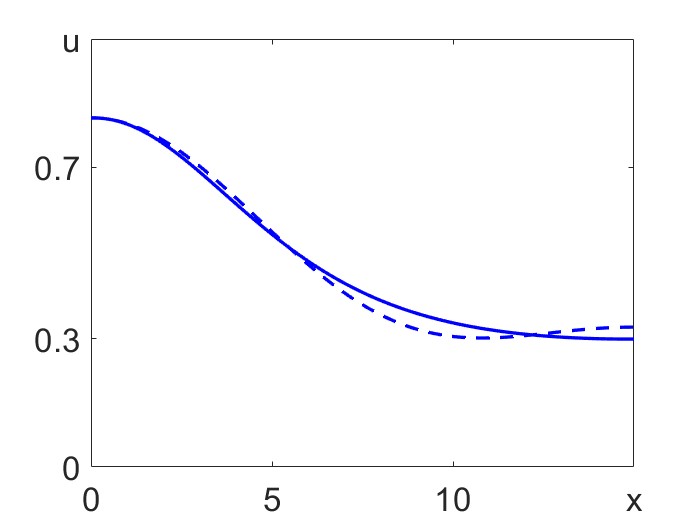}
\end{subfigure}
\begin{subfigure}{.29\textwidth}
\caption{}
\includegraphics[scale=.21]{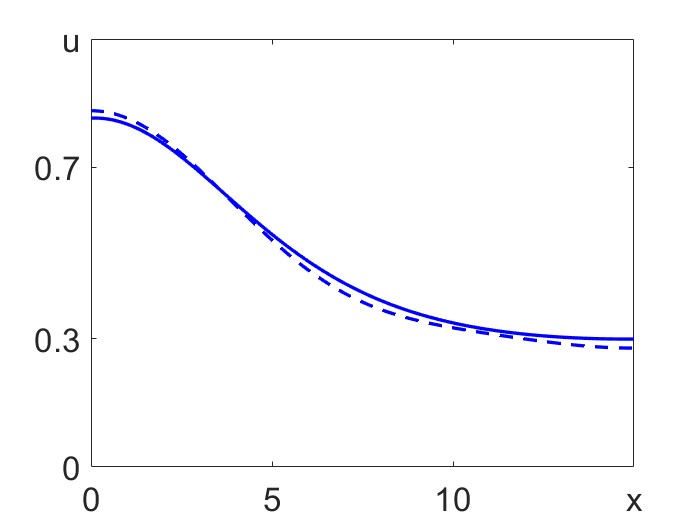}
\end{subfigure}
\begin{subfigure}{.27\textwidth}
	\caption{}
	\includegraphics[scale=.21]{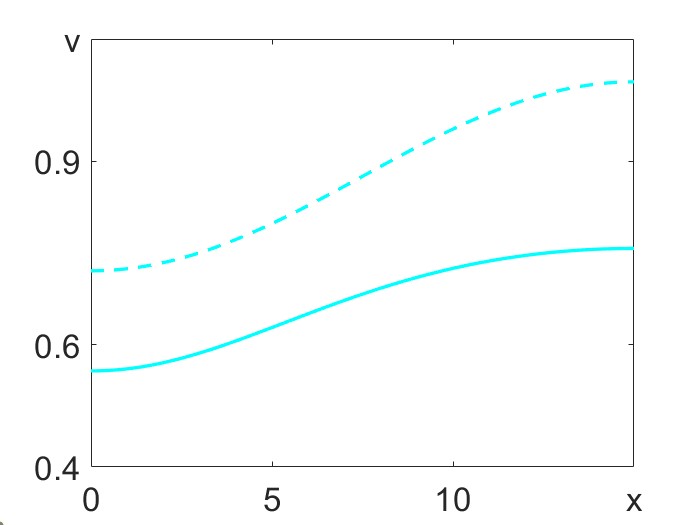}
\end{subfigure}
\hfill
\begin{subfigure}{.35\textwidth}
	\caption{}
	\includegraphics[scale=.21]{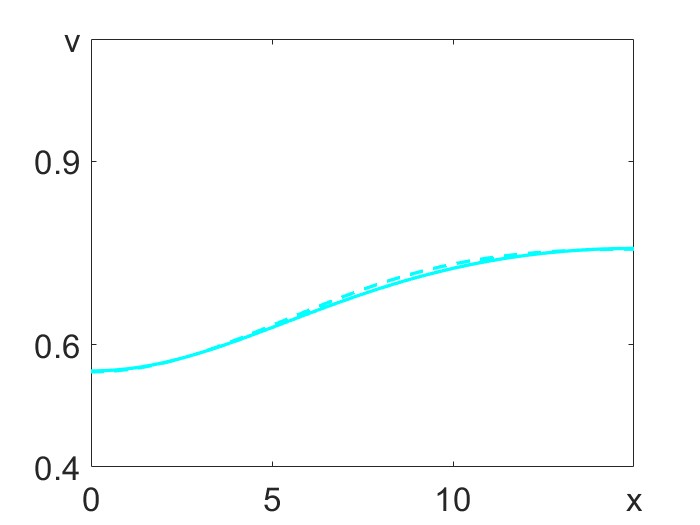}
\end{subfigure}
\begin{subfigure}{.29\textwidth}
	\caption{}
	\includegraphics[scale=.21]{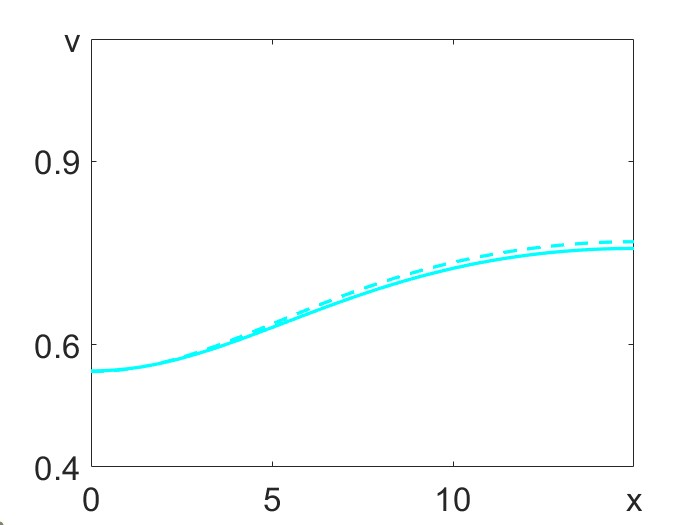}
\end{subfigure}
\caption{\em{\textbf{Comparison between the numerical and analytical profiles of $u(x)$ and $v(x)$.} \textbf{(a)}: Numerical $u$-profile (solid) versus analytical $u$-profile (dashed) obtained by truncating the system \eqref{Fsimeqs_twosp} at $M=1$. Error between the curves is $ER=1.3830$. \textbf{(b)}: Increasing the truncation to $M=2$ significantly reduces the error between the numerical and analytical $u$-profiles to $ER=0.0061$. \textbf{(c)}: Further increase in truncation to $M=3$ reduces the error between the two curves to $ER=0.0035$. Similar results for $v$-profiles: truncation at $M=1$ with error $ER=0.6566$ on panel \textbf{(d)}, $M=2$ with error $ER=0.0003$ on panel \textbf{(e)} and $M=3$ with error $ER=0.0009$ on panel \textbf{(f)}. Parameter values: $D_1=D_2=1, \mbox{ } r_1=r_2=0.1, \mbox{ } \chi=-10$ and $b_1=b_2=0.7$. }}
\label{twosp_a4prof}
\end{figure}

Fig.\ref{twosp_a4prof}{(a-c)} shows that solving the system ({\ref{Fsimeqs_twosp}) for $M=3$ is sufficient to accurately reproduce the $u(x)$ pattern, as increasing the truncation further to $M=4$ results in solving a system of 10 simultaneous equations, which is more computationally expensive and does not significantly affect the error between the two curves. This result matches the one obtained through the numerical integration of the profile from simulations, as panel $(a)$ in Figure \ref{ag_numerically} shows that only coefficients up to and including $\alpha_3$ influence the amplitude of the half spike. A similar analysis is performed for the pattern produced by the second species, $v(x)$, which has a smaller amplitude compared to that of the first species, $u(x)$ (see Fig.\ref{twosp_a4prof}{(d-f)}).


Patterns forming in the case when $b_2>1$ can also be found as solutions of the system \eqref{Fseriessystem}. Again for $M>0$ these solutions can't be expressed as an explicit function of model parameters. However, numerical  solutions can be found for fixed set of parameter values with the aid of Maple. We have made such calculations for the set of parameter values $D_1 = D_2 = 1$, $r_1 = r_2 = 0.1$, $\chi = -10$, $b_1 = 0.7$, $b_2 = 1.7$ and truncations $M=1, 2, 3$ and $4$ of the system \eqref{Fsimeqs_twosp}. In Table \ref{Maple_a4_100} obtained $\alpha$ coefficients (for a medium length $L = 10$, representing the half-wavelength of the numerically simulated pattern) are compared against coefficients found by Fourier decomposition of the numerical solution represented in Fig.\ref{ag_numerically}$(b)$.
It is evident that the truncation at $M=1$ gives a significant discrepancy between the numerical and analytical coefficients, indicating that the truncation needs to be increased to include contributions from higher frequency modes. The coefficients obtained by truncating the series at $M=2, 3$ and $4$ provide much better approximations of the numerical ones, suggesting that this analytical solution is a closer representation of the numerical solution. It is evident from comparison of the $u$-profiles, representing  solutions of \eqref{Fsimeqs_twosp} with the profile obtained in numerical simulations. The discrepancy between these profiles can be calculated using equation \eqref{erfor}.
\begin{table}[h]
	\begin{center}
		\begin{tabular}{| l | l | l | l | l | l | } 
			\hline
			$L=10$ & $\alpha_0$  & $\alpha_1$ & $\alpha_2$ & $\alpha_3$ & $\alpha_4$  \\[7pt] \hline
			Numerical & 0.3234 & -0.4206 & 0.2015 & -0.0779 & 0.0283   \\[7pt] \hline
			$M=1$ & 0.6419 & -0.5806 & 0 & 0 & 0 \\[7pt] \hline
			$M=2$ & 0.3248 & -0.4128 & 0.2878 & 0 & 0 \\[7pt] \hline
			$M=3$ & 0.3139 & -0.4367 & 0.2208 & -0.1080 & 0 \\[7pt] \hline
			$M=4$ & 0.3208 & -0.4419 & 0.2241 & -0.0916 & 0.0342 \\[7pt] \hline
		\end{tabular}
	\end{center}
	\caption{\em{\textbf{First five Fourier coefficients describing $u$-profile in the case when the competition is not weak.} Numerical coefficients are obtained by Fourier decomposition of the simulated profile $u(x)$, while analytical coefficients are obtained by solving system (\ref{Fsimeqs_twosp}) truncated at $M=1, 2, 3$ and $4$. Parameter values: $D_1=D_2=1, \mbox{ } r_1=r_2=0.1, \mbox{ } \chi=-10, \mbox{ }b_1=0.7, \mbox{ }b_2=1.7$ and medium length $L=10$.}}
	\label{Maple_a4_100}
\end{table}

Figure \ref{twosp_a4prof_100} provides a graphical comparison between the numerical profile from simulations and the analytical profiles given as solutions of the system \eqref{Fsimeqs_twosp}. As expected, panel $(a)$ shows a significant difference between the profile from simulations and the profile given by the Fourier series truncated at $M=1$, with the discrepancy between the two curves being $ER=1.3928$. In panel $(b)$, with $M=2$ such that coefficient $\alpha_2 \neq 0$, there is a considerable improvement in the shape of the analytical profile, as well as a reduction in the discrepancy between the two curves, which is now $ER=0.0717$. The analytical profile in panel $(c)$ represents the Fourier solution when the series for $u$ is truncated at $M=3$, and the error between the curves is further reduced to $ER=0.0155$. Increasing the truncation further in panel $(d)$ does not significantly improve the error between the two curves, $ER=0.0113$, but it does produce a smoother analytical profile. Since increasing to $M=4$ did not have a substantial effect on the discrepancy between the curves, we conclude that truncation at $M=4$ is sufficient to capture the analytical solution corresponding to species $u(x)$ that emerges when $b_2>1$.


\begin{figure}[h]
	\centering
	\begin{subfigure}{.16\textwidth}
		\caption{}
		\includegraphics[scale=.15]{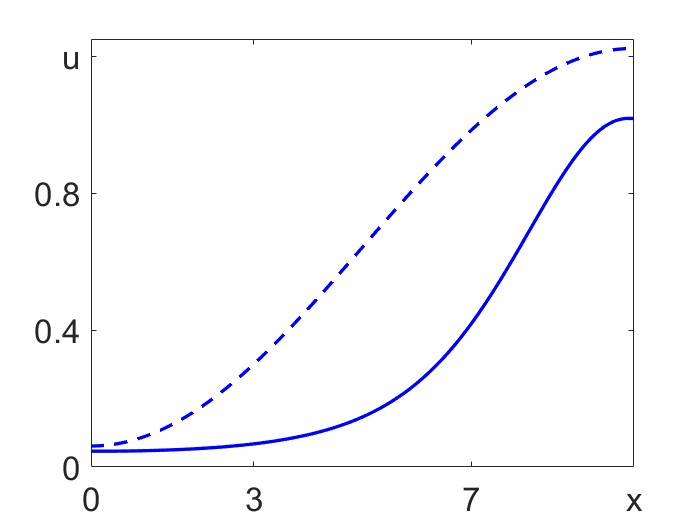}
	\end{subfigure}
	\hfill
	\begin{subfigure}{.25\textwidth}
		\caption{}
		\includegraphics[scale=.15]{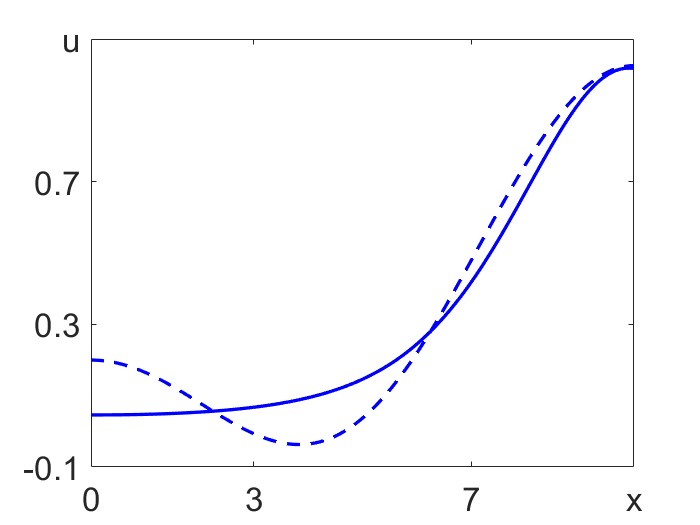}
	\end{subfigure}
	\begin{subfigure}{.16\textwidth}
		\caption{}
		\includegraphics[scale=.15]{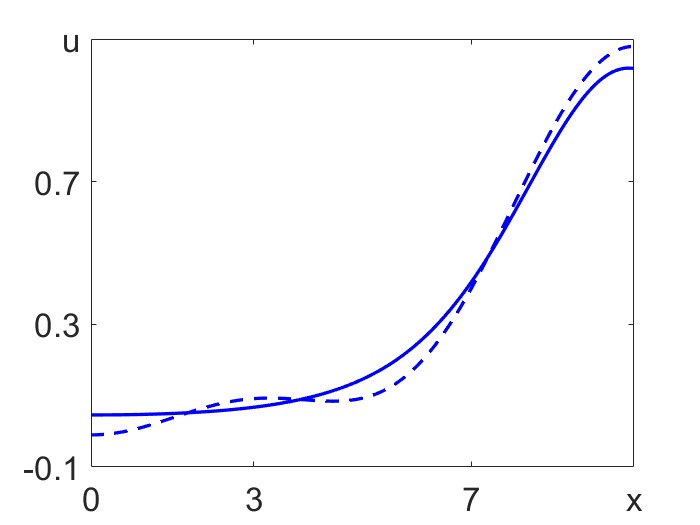}
	\end{subfigure}
	\hfill
	\begin{subfigure}{.25\textwidth}
		\caption{}
		\includegraphics[scale=.15]{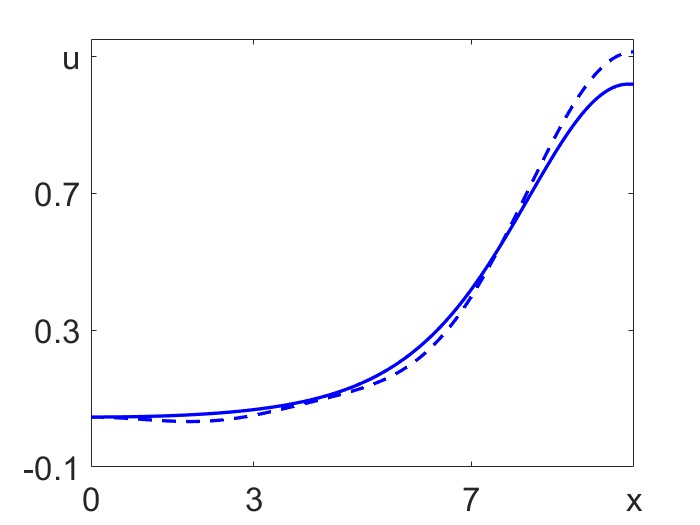}
	\end{subfigure}
	\begin{subfigure}{.16\textwidth}
	\caption{}
	\includegraphics[scale=.15]{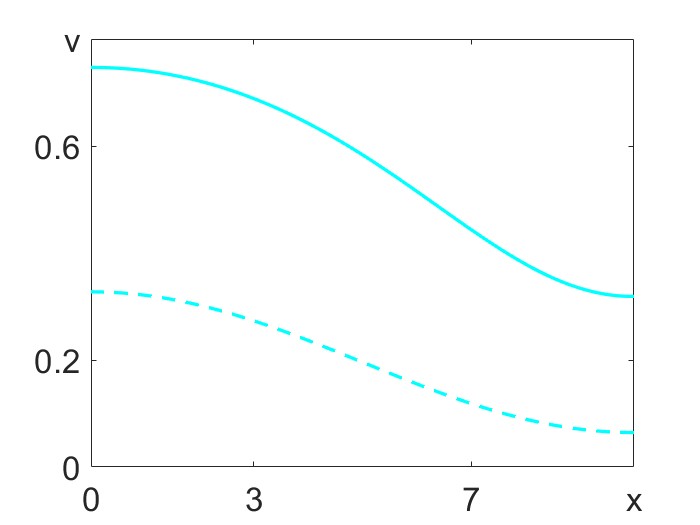}
\end{subfigure}
\hfill
\begin{subfigure}{.25\textwidth}
	\caption{}
	\includegraphics[scale=.15]{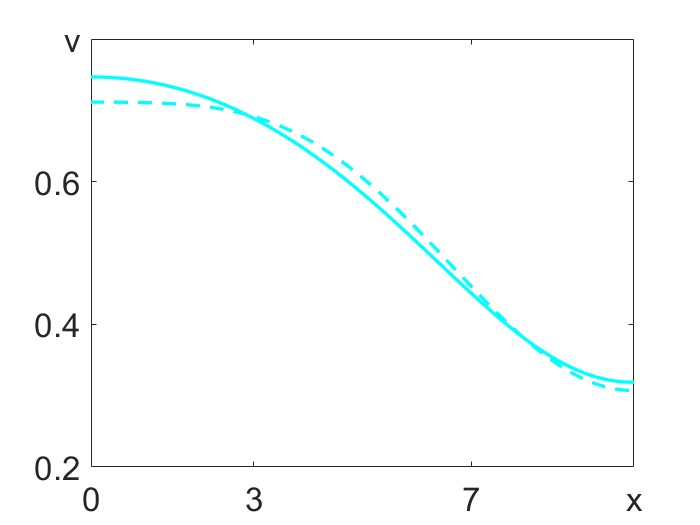}
\end{subfigure}
\begin{subfigure}{.16\textwidth}
	\caption{}
	\includegraphics[scale=.15]{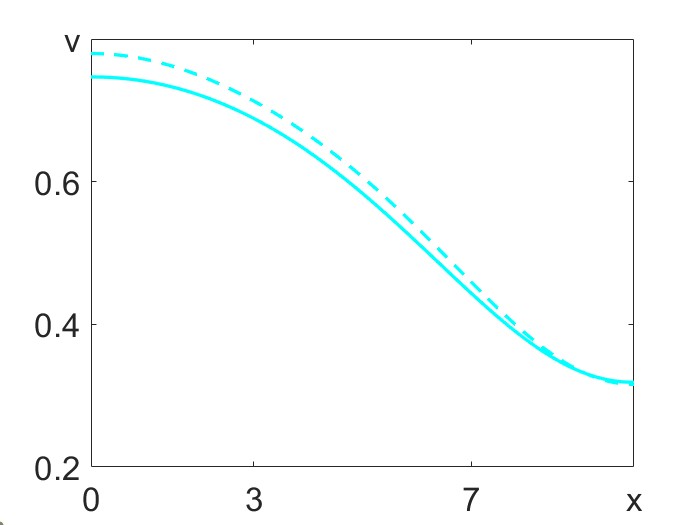}
\end{subfigure}
\hfill
\begin{subfigure}{.25\textwidth}
	\caption{}
	\includegraphics[scale=.15]{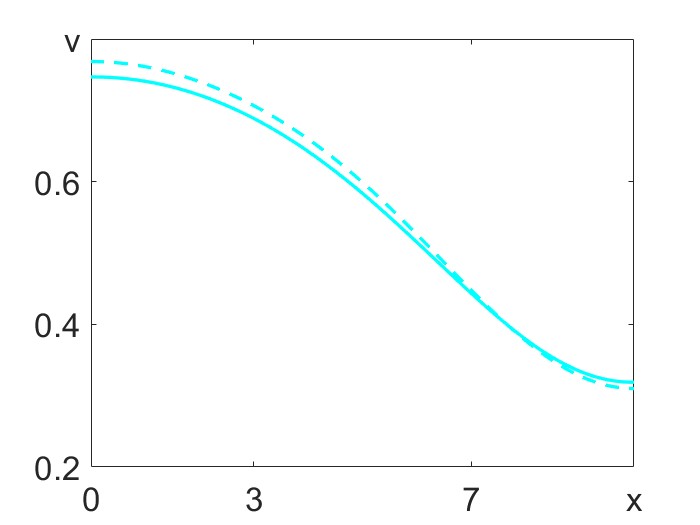}
\end{subfigure}
	\caption{\em{\textbf{Comparison between the numerical and analytical profiles of $u(x)$ and v(x).} \textbf{(a)}: Numerical $u$-profile (solid) versus analytical $u$-profile (dashed) obtained by truncating the system \eqref{Fsimeqs_twosp} at $M=1$. Error between the curves is $ER=1.3928$. \textbf{(b)}: Increasing the truncation to $M=2$ significantly reduces the error between the numerical and analytical curves to $ER=0.0717$. \textbf{(c)}: Further increase in truncation to $M=3$ reduces the error between the two curves further to $ER=0.0155$. \textbf{(d)}: Increasing the truncation to $M=4$ does not have a significant effect on the difference between the two curves, $ER=0.0113$, but makes a smoother analytical profile. Fixed parameters: $D_1=D_2=1, \mbox{ } r_1=r_2=0.1, \mbox{ } \chi=-10$ and $b_1=0.7, \mbox{ }b_2=1.7$. 
Similar results for $v$-profiles: truncation at $M=1$ with error $ER=1.3482$ on panel \textbf{(e)}, $M=2$ with error $ER=0.0031$ on panel \textbf{(f)}, $M=3$ with error $ER=0.0047$ on panel \textbf{(g)} and $M=4$ with error $ER=0.0019$ on panel \textbf{(h)}. Parameter values: $D_1=D_2=1, \mbox{ } r_1=r_2=0.1, \mbox{ } \chi=-10$ and $b_1=b_2=0.7$.}}
	\label{twosp_a4prof_100}
\end{figure}

This section has demonstrated that nonlinear Fourier analysis is a powerful method for describing the periodic patterns forming both in weak and non-weak competition cases. We obtained Fourier series coefficients for the two species, $u$ and $v$, by first numerically investigating the pattern obtained from simulations and then solving a set of simultaneous equations to find the same coefficients analytically. A comparison between numerical and analytical results has shown that to accurately describe the pattern containing half spike a truncation of $M=3$ is necessary in the weak competition case and $M=4$ in the case when $b_2>1$. This result was demonstrated for a certain fixed set of model parameters. Next task is to investigate the impact of model parameters on amplitude and wavelength of forming periodic patterns. 

\section{Effect of model parameters on amplitude and wavelength of periodic patterns}
 
In this section, the effect of model parameters on amplitude and wavelength is investigated. For simplicity in the analytical Fourier analysis, the medium length, $L$, is constrained to focus on the characteristics of half-spike patterns; however, the same results hold for patterns with multiple spikes. The effect of parameters is explored by running simulations for different medium lengths with a step size of 1, integrating the profile according to (\ref{twospecies_Fouriersol}), and observing how the fastest-growing nodes, $\alpha_1$ and $\gamma_1$, change with $L$. The lengths that produce the largest node correspond to the most unstable wavelength, which is the fastest-growing medium length, resulting in the highest pattern amplitude. Numerical results have been verified using analytical Fourier analysis by solving system (\ref{Fsimeqs_twosp}) for $M=4$ and similarly investigating the effects of $L$ on $\alpha_1$ and $\gamma_1$, corresponding to the profiles  $u(x)$ and $v(x)$, respectively.

Diffusion of species plays an important role in spatial distribution and, consequently, pattern formation. Therefore, the effects of $D_1$ and $D_2$, corresponding to species  $u$ and $v$, respectively, are investigated first. Results for how the most unstable wavelength and the fastest-growing mode change for different diffusion coefficients are represented graphically in Figure \ref{D1D2_nvag}.

\begin{figure}[h]
\centering
\begin{subfigure}{.16\textwidth}
\caption{}
\includegraphics[scale=.15]{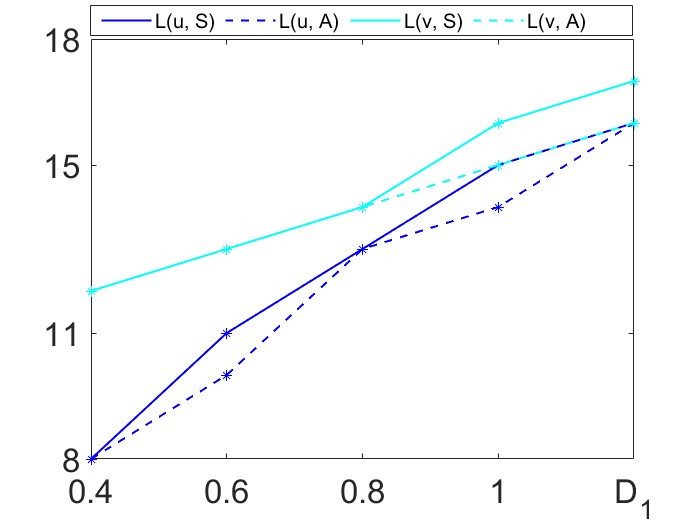}
\end{subfigure}
\hfill
\begin{subfigure}{.25\textwidth}
\caption{}
\includegraphics[scale=.15]{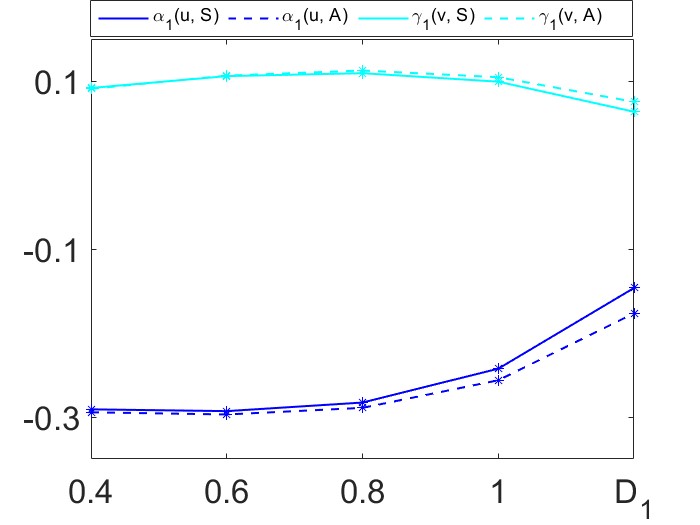}
\end{subfigure}
\begin{subfigure}{.16\textwidth}
\caption{}
\includegraphics[scale=.15]{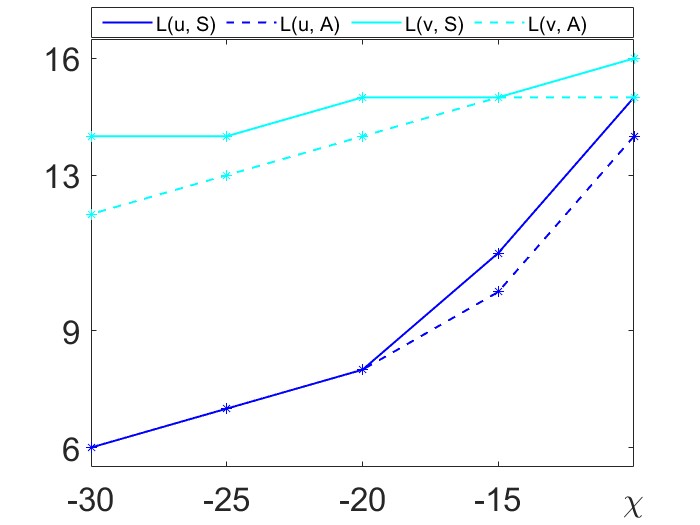}
\end{subfigure}
\hfill
\begin{subfigure}{.25\textwidth}
\caption{}
\includegraphics[scale=.15]{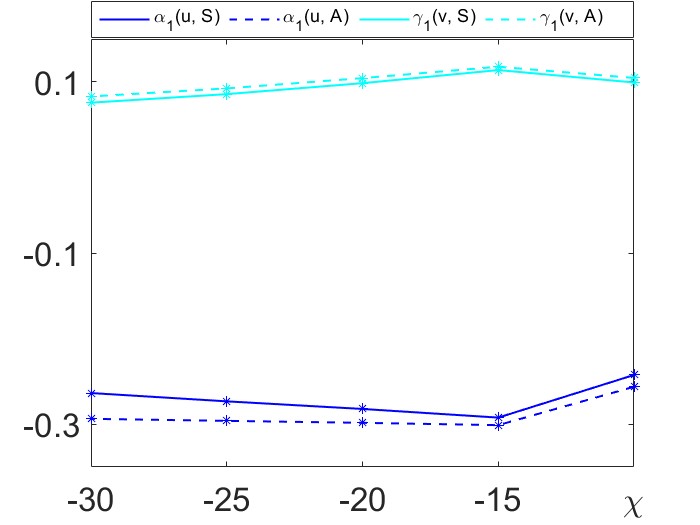}
\end{subfigure}
\begin{subfigure}{.16\textwidth}
	\caption{}
	\includegraphics[scale=.15]{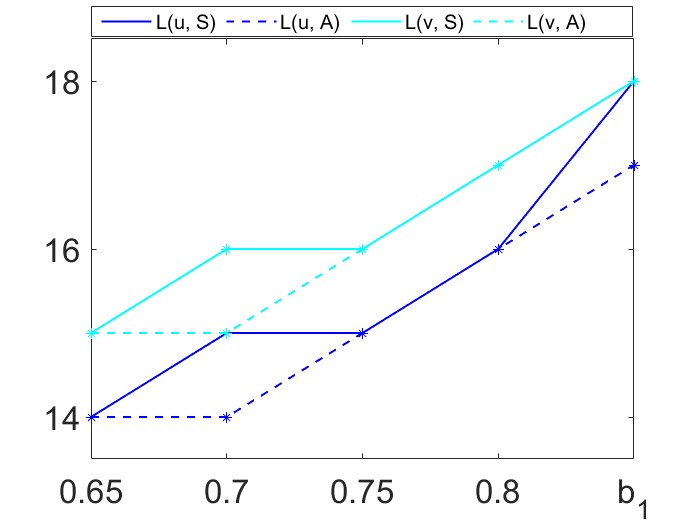}
\end{subfigure}
\hfill
\begin{subfigure}{.25\textwidth}
	\caption{}
	\includegraphics[scale=.15]{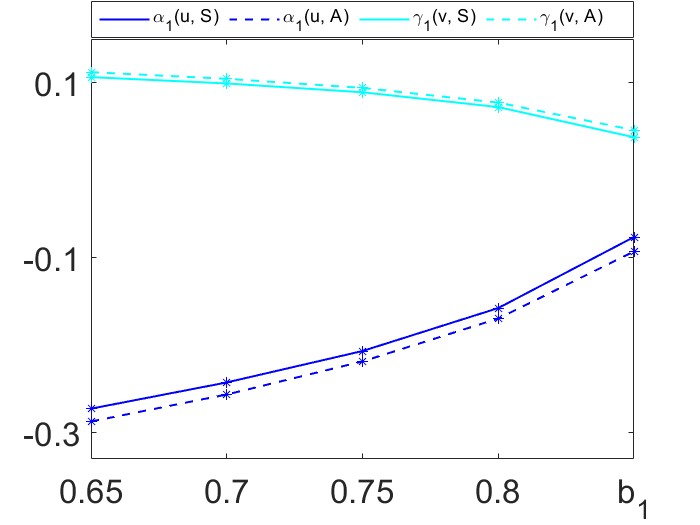}
\end{subfigure}
\begin{subfigure}{.16\textwidth}
	\caption{}
	\includegraphics[scale=.15]{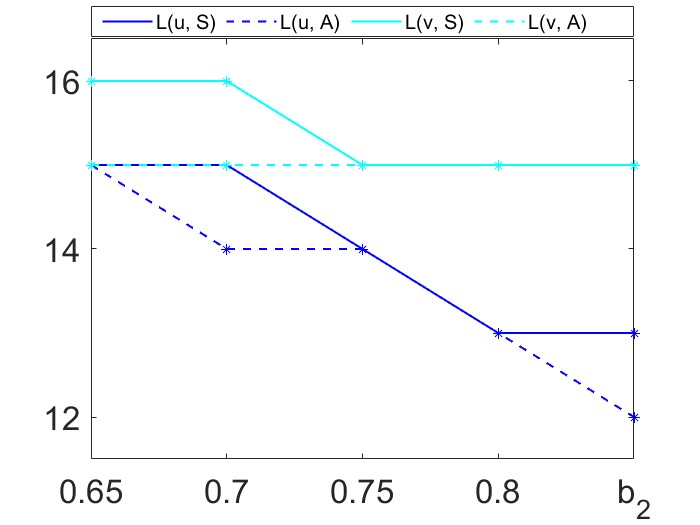}
\end{subfigure}
\hfill
\begin{subfigure}{.25\textwidth}
	\caption{}
	\includegraphics[scale=.15]{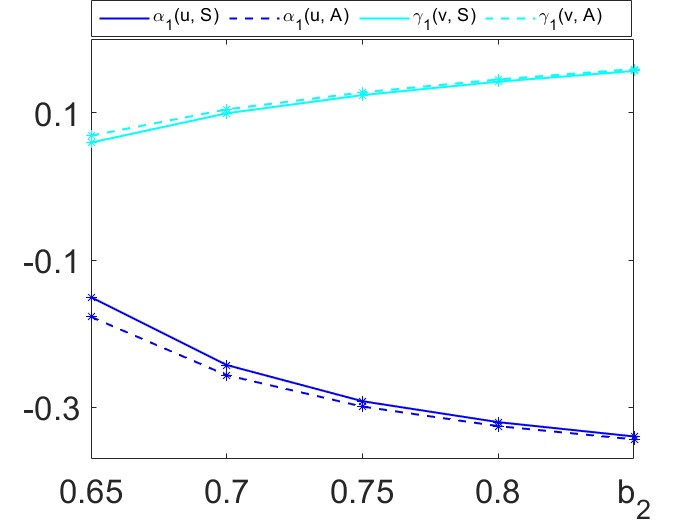}
\end{subfigure}
\begin{subfigure}{.16\textwidth}
	\caption{}
	\includegraphics[scale=.15]{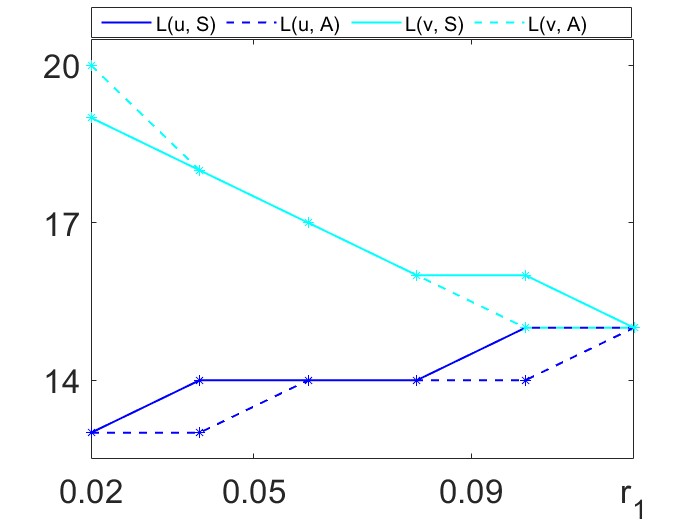}
\end{subfigure}
\hfill
\begin{subfigure}{.25\textwidth}
	\caption{}
	\includegraphics[scale=.15]{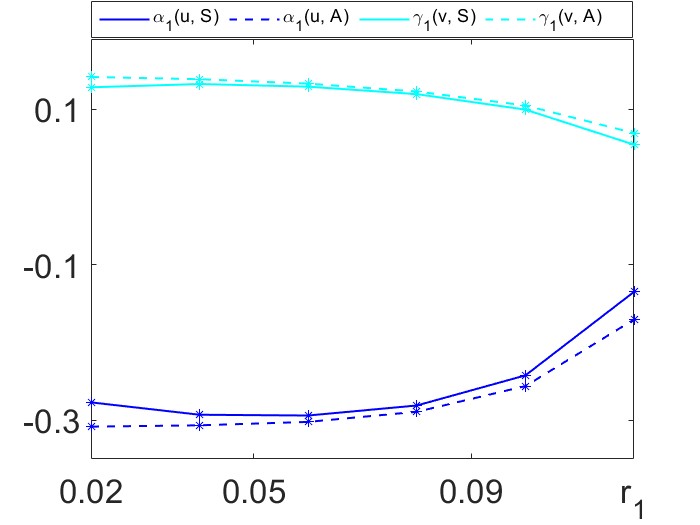}
\end{subfigure}
\begin{subfigure}{.16\textwidth}
	\caption{}
	\includegraphics[scale=.15]{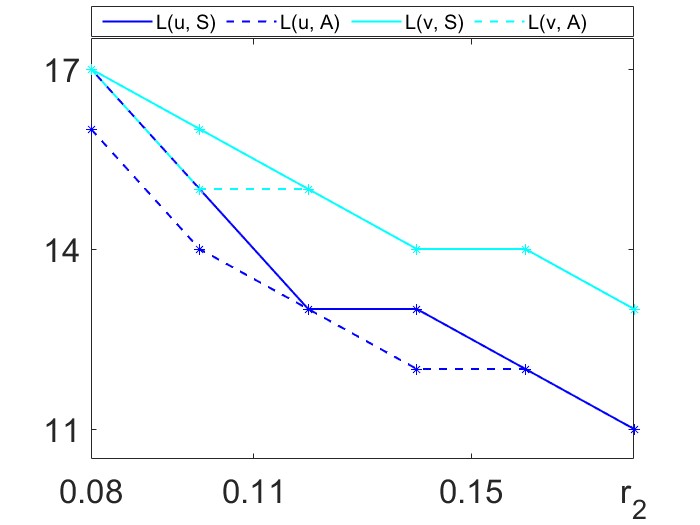}
\end{subfigure}
\hfill
\begin{subfigure}{.25\textwidth}
	\caption{}
	\includegraphics[scale=.15]{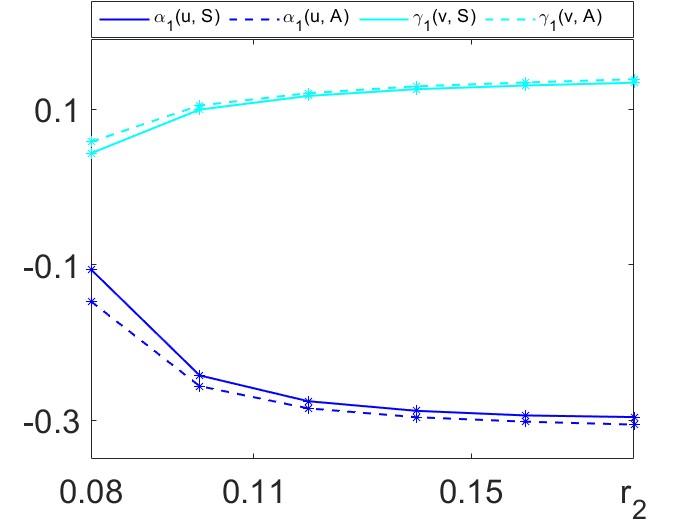}
\end{subfigure}
\caption{\em{\textbf{Dependence of the wavelength and amplitude of periodic patterns on model parameters.}  The effect of \textbf{(a-b)} the diffusion coefficient, $D_1$; \textbf{(c-d)} chemotaxtic sensitivity, $\chi$; \textbf{(e-f)} competition rate, $b_1$; \textbf{(g-h)} competition rate, $b_2$; \textbf{(i-j)} proliferation rate, $r_1$, and \textbf{(k-l)} proliferation rate, $r_2$. Default set of parameters: $D_1=D_2=1$, $r_1=r_2=0.1$,  $b_1=b_2=0.7$ and $\chi=-10$. 
}}
\label{D1D2_nvag}
\end{figure}

Figure \ref{D1D2_nvag} shows the effects of diffusion on the amplitude and wavelength of the half-spike pattern formed by the two species  $u(x)$ and $v(x)$. $L(u, S)$ refers to the most unstable wavelength corresponding to species $u(x)$, obtained by running numerical simulations for different medium lengths $L$ and determining which length yields the largest $\alpha_1$. Similarly, $L(u, A)$ represents the most unstable wavelength obtained analytically by solving system (\ref{Fsimeqs_twosp}) truncated at $M=4$. $\alpha_1(n, S)$ denotes the largest, fastest-growing node obtained from simulations, while $\alpha_1(n, A)$ refers to the fastest-growing node obtained analytically. If the absolute value of $\alpha_1$ increases, the pattern's amplitude also increases, and if the absolute value of $\alpha_1$ decreases, the amplitude decreases accordingly. The same notations apply for species $v(x)$. Panels $(a)$ and $(b)$ demonstrate that the most unstable wavelength increases with increasing diffusion of  $u$, represented by $D_1$, while the amplitude of the pattern decreases as $D_1$ increases. An interesting characteristic of the patterns corresponding to species  $u$ and $v$ is that, for small values of $D_1$, the most unstable wavelengths of the two species differ significantly. This suggests that irregularities in the patterns are expected if the diffusion of  $u$ is small. Similarly, panels $(c)$ and $(d)$ show that an increase in $D_2$ leads to an increase in the most unstable wavelength, while the amplitude of the pattern decreases. Additionally, irregularities are expected to appear for small values of $D_2$. All panels display good correlation between numerical and analytical results, indicating that accurate solutions to system (\ref{twospecieschem}) can be found by solving (\ref{Fsimeqs_twosp}).

Pattern formation in a two-species system arises as a result of chemorepulsion. As a reminder, chemorepulsion needs to be strong enough for system (\ref{twospecieschem}) to become unstable under perturbation. Therefore, the effects of chemotaxis on wavelength and amplitude are investigated. Numerical and analytical results are represented graphically in panels $(c)$ and $(d)$ of Figure \ref{D1D2_nvag}, which illustrates the effects of chemorepulsion on the amplitude and wavelength of patterns resulting from the disturbance of the coexistence steady state of system (\ref{twospecieschem}). Panel $(c)$ shows that as chemorepulsion increases, the most unstable wavelength decreases. Additionally, with increasing chemorepulsion, a significant difference arises between the most unstable wavelength of species  $u$ and that of species $v$, indicating that strong chemorepulsion leads to irregularities in pattern formation. Panel $(d)$ demonstrates that the maximal amplitude of the pattern formed by both species occurs at $\chi=-15$. As chemorepulsion deviates from this value, the amplitude of the pattern decreases. It is worth noting that the analytical Fourier solutions predict a slightly larger pattern amplitude than the one obtained from simulations, although there is a very good approximation between the numerical and analytical results.

Next, the effects of interspecific competition are examined, as these parameters play a crucial role in the formation of periodic patterns around the coexistence steady state. Figure \ref{D1D2_nvag} shows the effects of interspecific competition on the wavelength and amplitude of the half-spike pattern formed by an infinitesimal perturbation of the coexistence steady state. In panels $(e)$ and $(f)$, the effects of the competition of $v$ on  $u$ are investigated. As $b_1$ increases, the most unstable wavelength also increases, while the amplitude of the pattern decreases. Conversely, panels $(g)$ and $(h)$ show the effects of the competition of  $u$ on $v$, and in this case, as $b_2$ increases, the most unstable wavelength decreases, and the amplitude of the pattern increases. Additionally, in both cases, there are differences between the most unstable wavelength of  $u$ and that of $v$, which suggests that pattern irregularities are likely to appear. It is also evident that $b_1$ and $b_2$ have opposite effects on the characteristics of the pattern: $b_2$ enhances pattern formation, leading to larger pattern amplitudes, while $b_1$ leads to smaller ones.

The final set of parameters to investigate concerns the effects of reproduction on pattern formation. As shown in Figure \ref{TP_b}, periodic patterns are more likely to appear when species  $u$ has a smaller reproduction rate, $r_1$, and species $v$ has a larger reproduction rate, $r_2$. This leads to the prediction that an increase in $r_1$ would decrease the amplitude of the pattern, while an increase in $r_2$ would enhance it. Figure \ref{D1D2_nvag} illustrates the effects of species reproduction on pattern characteristics. Panels $(i)$ and $(j)$ show that as the reproduction rate of the first species increases, the most unstable wavelength decreases, along with the amplitude of the pattern. These results are consistent with previous findings, which indicated that, for this set of fixed parameters, pattern formation is not possible for $r_1 > 0.15$. Similarly, panels $(k)$ and $(l)$ show that as $r_2$ increases, the most unstable wavelength also decreases; however, in this case, the amplitude of the pattern increases. This aligns with earlier results showing that, for this set of fixed parameters, pattern formation is not possible for small values of $r_2 < 0.8$. See Figure \ref{TP_b}.

This section has focused on understanding how patterns are affected by changes in parameters, specifically examining how characteristics such as amplitude and wavelength are influenced. This has been demonstrated through methods of numerical and analytical Fourier analysis, which show a strong correlation between the two sets of results. This indicates that truncating system (\ref{Fsimeqs_twosp}) at $M=3$ and solving analytically accurately captures the characteristics of a half-spike pattern. Additionally, it has been shown that the parameters which enhance the amplitude of the pattern as they increase are $b_2$ and $r_2$, while all other parameters, $D_1$, $D_2$, $\chi$, $b_1$, and $r_1$, reduce the amplitude as they increase. It has also been observed that species  $u$ and $v$ exhibit different wavelengths under certain parameter regimes, suggesting that patterns formed from the coexistence steady state of model (\ref{twospecieschem}) are likely to display irregularities.

This concludes the analysis of pattern formation around the coexistence steady state, and the next section focuses on pattern formation around the extinction steady states.

\section{Discussion}

Understanding interactions between two or more bacterial species is crucial in the early stages of colonisation and biofilm formation. Systems of two competing bacterial species are often referred to as Lotka-Volterra models, and the stability of the steady states has been thoroughly analysed \cite{Lotka, Volterra}. In this chapter, the Lotka-Volterra model and the Keller-Segel model \cite{Ks, Keller} have been combined to investigate a system of two bacterial species,  $u$ and $v$, where the latter produces a chemical agent, $c$, which chemotactically affects the former. These interactions are modelled using a system of three partial differential equations (\ref{twospecieschem}), which has four steady states:
\[\displaystyle (u,\mbox{ }v, \mbox{ }c)=\left\{(0, 0, 0),\mbox{ } (1, 0, 0), \mbox{ }(0, 1, 1), \mbox{ }\left(\frac{b_1-1}{b_1b_2-1},\mbox{ } \frac{b_2-1}{b_1b_2-1}, \mbox{ } \frac{b_2-1}{b_1b_2-1}\right)\right\}, \]
where the trivial steady state is always unstable, while the other three steady states are stable in the well-mixed system under specific conditions. Classical Turing pattern analysis states that pattern formation occurs if an initially stable steady state in the well-mixed system is driven unstable by a perturbation in the full reaction-diffusion-advection system \cite{Turing}. In this chapter, methods of linear analysis and nonlinear Fourier analysis have been used to demonstrate the formation of periodic patterns around both the coexistence steady state $\displaystyle \left(u, v, c\right)=\left(\frac{b_1-1}{b_1b_2-1}, \frac{b_2-1}{b_1b_2-1}, \frac{b_2-1}{b_1b_2-1}\right)$ and the extinction steady state $(u, v, c)=(1, 0, 0)$, as a result of aggregation initiated by a breakdown of stability.

As demonstrated in the introductory chapter, the coexistence steady state is stable in the well-mixed system if the interspecific competition between the species is less than 1, i.e., $b_1, b_2 < 1$. This identifies the region of interest for the formation of stationary periodic patterns. The first section of this chapter focused on showing that model (\ref{twospecieschem}) at coexistence can be driven unstable by infinitesimal perturbations within this region. Since finding analytical conditions for instability proved difficult, the Routh-Hurwitz criteria \cite{Hurwitz} were used with fixed parameters to identify a domain, $\delta_2$, in which the system is unstable and exhibits pattern formation. Using numerical simulations, we have shown that for competition rates $b_1, b_2 \in \delta_2$, stationary periodic patterns arise from infinitesimal disturbances as a result of chemorepulsion. An important question that followed was how to ensure we consider the largest domain of instability, $\delta_2$, in order to capture all competition rates that lead to a breakdown of stability and initiate aggregation. To address this, we found the most unstable wavenumber, $k$, which represents the fastest-growing mode. By examining the intersection point between $a_3=0$ and $(a_3)_k=0$, we determined the most unstable wavenumber $k$, as well as the minimum competition rates $b=b_1=b_2$ that lead to pattern formation. This prompted further analysis of how the minimum $b$ that leads to instability changes with other parameters. This analysis was conducted through numerical simulations and analytics, and the results showed strong correlation. As diffusion increases, the minimum requirement for $b$ also increases, meaning that faster diffusing species require stronger competition to initiate aggregation. As the strength of chemorepulsion increases, the minimum requirement for $b$ decreases. Additionally, an increase in the reproduction rate, $r_1$, of the repelled species  $u$ leads to an increase in the minimum requirement for $b$, whereas an increase in the reproduction rate, $r_2$, of the species producing the repellent reduces the minimum requirement for $b$. This means that as species $v$ reproduces faster, more chemotactic agent is produced, and aggregation is initiated even with weaker competition between species.

Classical Turing pattern analysis has provided conditions under which the system is driven unstable by perturbation, leading to pattern formation. However, it offers no information about the characteristics of the pattern, such as amplitude and wavelength. Using nonlinear Fourier pattern analysis, the aim is to represent solutions of model (\ref{twospecieschem}) as Fourier series and identify the fastest-growing mode. Initially, the analysis was carried out numerically by integrating a simulated profile to find the coefficients corresponding to the Fourier series solution. It has been shown that, in order to accurately reproduce a pattern formed by species  $u$ containing $i$ half-spikes, $\alpha_i$, $\alpha_{2i}$, and $\alpha_{3i}$ must be determined, and similarly for patterns produced by $v$ and $c$. For a pattern with many spikes, these coefficients are easy to obtain numerically, but much more difficult to derive analytically. Therefore, to obtain the analytical coefficients corresponding to the Fourier series solution, the medium length was reduced to obtain patterns containing at most one half-spike or a full spike. Analytical Fourier series are obtained by truncating at some $M$, which represents the highest frequency mode, and the goal is to determine which $M$ provides a good approximation of the profile obtained from simulations. It has been shown that truncating at $M=1$ results in a large discrepancy between the numerical and analytical profiles. However, including one more term in the series and truncating at $M=2$ produces a much better approximation. Results were verified for truncations up to $M=4$, and we concluded that this is an acceptable truncation parameter for representing a half-spike pattern, based on the discrepancy between numerical and analytical solutions. The effects of model parameters on the most unstable mode have also been investigated using numerical simulations and Fourier analysis. It has been shown that the analytically obtained Fourier solution for a half-spike pattern accurately captures the characteristics of the numerical solution. Moreover, for certain parameter regimes, the two species  $u$ and $v$ have different most unstable wavelengths, a characteristic that explains why pattern irregularities appear in a system of two bacterial species. The most unstable wavelength is the one corresponding to the fastest-growing mode. This has been used to show that certain model parameters, such as $b_2$ and $r_2$, enhance pattern formation as they increase, meaning that the amplitude of the pattern increases. On the other hand, most other model parameters reduce the amplitude of the pattern as they increase.

The second part of the chapter has focused on pattern formation from the extinction steady state $(u, v, c) = (1, 0, 0)$ due to finite amplitude disturbance. This extinction steady state is stable in the well-mixed system if $b_2 > 1$ and remains stable under the same condition in the full reaction-diffusion-advection system. According to linear analysis, stationary periodic patterns do not emerge if this steady state is perturbed by an infinitesimal disturbance. However, numerical simulations have shown that if species $v$ is perturbed by a finite amplitude disturbance, $\tilde{v}$, stationary periodic patterns emerge, provided the disturbance is large enough. Initially, this phenomenon was analysed using computational simulations, and one of the first questions addressed was how the minimum perturbation amplitude $\tilde{v}$ is affected by parameters. It was shown that if diffusion increases, the minimum perturbation amplitude also increases. This means that for more diffusive species, a larger density of species $v$ needs to be introduced into the system to prevent it from decaying back to a homogeneous state before aggregation is initiated. As expected, an increase in chemorepulsion strength decreases the minimum perturbation amplitude required to initiate aggregation. For an increase in $b_1$, which represents the competition of $v$ on  $u$, there is a decrease in $\tilde{v}$, since this indicates that $v$ is more competitive for resources and is more likely to survive and produce enough chemical agent to initiate aggregation. Conversely, for an increase in $b_2$, which represents the competition of  $u$ on $v$, there is an increase in the minimum perturbation amplitude required to initiate aggregation. Similarly, an increase in the reproduction rate of the first species, $r_1$, results in an increase in $\tilde{v}$, whereas an increase in the reproduction rate of the second species, $r_2$, decreases the minimum perturbation amplitude required to initiate aggregation. Solutions around this steady state have been investigated using nonlinear Fourier analysis, following the same procedure as for patterns emerging from the coexistence steady state. In this case, exact analytical solutions were found by solving a system of simultaneous equations (\ref{Fsimeqs_twosp}) for fixed parameters, and it was shown that truncating the series at $M=4$ produces accurate results with small discrepancies between numerical and analytical profiles. The most important result of this section is that, contrary to linear analysis, stationary periodic patterns can emerge from the steady state $(n, v, c) = (1, 0, 0)$ due to finite amplitude disturbances.

This chapter has focused on the analysis of stationary periodic patterns around the coexistence steady state and one of the extinction steady states, $(u, v, c) = (1, 0, 0)$. However, there is one more extinction steady state, $(u, v, c) = (0, 1, 1)$, which is stable in the well-mixed system if $b_1 > 1$. Classical Turing pattern analysis suggests that stationary periodic patterns can emerge from this steady state if the system can be driven unstable by perturbation in the full reaction-diffusion-advection system under the same condition. However, according to linear analysis, this steady state remains stable when perturbed in the full diffusive system for $b_1 > 1$, and pattern formation is not possible. This result is consistent with findings from computational simulations and nonlinear Fourier analysis, even in the presence of finite amplitude disturbances.

In conclusion, this chapter has analysed pattern formation in a system consisting of two bacterial species interacting with a chemical agent. Pattern formation has been proven to emerge from two steady states, and characteristics of these patterns have been identified using nonlinear Fourier analysis. Analytical results have been supported by numerical computations.

\end{document}